\documentclass[12pt]{article}
\usepackage{pifont}

\usepackage{amsfonts}
\usepackage{latexsym}
\usepackage{amsmath}
\usepackage{amssymb}
\usepackage{color}

\def\z{\zeta}
\def\a{\alpha}
\def\ep{\varepsilon}

\def\d{\delta}

\def\ph{\phi}
\def\g{\gamma}
\def\l{\lambda}

\def\b{\beta}

\title{Backward stochastic viability property with jumps and applications to the comparison theorem for multidimensional BSDEs with jumps }

\date{June 8, 2010}

\author{ Xuehong Zhu}

\begin{document}

\maketitle

\begin{abstract}
In this paper, we study conditions under which the solutions of a
backward stochastic differential equation with jump remains in a
given set of constrains.
This property is the so-called "viability property". As an application, we study the comparison theorem for multidimensional BSDEs with jumps. \\
\par  $\textit{Keywords: Viability property; BSDE; Comparison theorem.}$
\end{abstract}

%% \linenumbers

%% main text

\section{Introduction}\label{sec:intro}
\qquad Viability properties for stochastic differential  equations
(SDEs) and inclusions had been introduced and studied by Aubin-Da
Prato in [2], [3] and [4]. The key point of their work was a
"stochastic tangent cone" which generalized the well-known
Bouligand's contingent cone used for deterministic systems.
Buckdahn-Peng-Quincampoix-Rainer [6] used a new method to get the
necessary and sufficient condition for the viability property of
SDEs with control. They related viability with a kind of optimal
control problem and applied the comparison theorem for viscosity
solutions to some H-J-B equation. In 2000, renewed interests arouse
in the study of viability for backward stochastic differential
equations (BSDEs) in the paper of Buckdahn-Quincampoix-
R\v{a}\c{s}canu [7]. Their approach differed from the above
mentioned methods and based on the convexity of the distance
function of $K$ (when $K$ is a closed convex set). This enabled them
to deduce some condition in differential form on the distance
function of $K$ which is necessary as well as sufficient.

In the present paper, we shall adopt the approach of Buckdahn et al.
[7] to study the backward stochastic viability property (BSVP) with
jumps.

As we know, the comparison theorem and the converse comparison
theorem for real-valued BSDEs with or without jumps had been studied
by many mathematicians. We refer to Peng [13], Pardoux-Peng [12], El
Karoui-Peng-Quenez [9], Coquet-Hu-Memin-Peng [8],
Barles-Buckdahn-Pardoux [5], Wu [15], Royer [14] and their
references. In 2006, Hu-Peng [10] first use the BSVP to study the
comparison theorem for multidimensional BSDEs and get necessary and
sufficient condition.

So We can apply our results of the BSVP with jumps to study the
comparison theorem for multidimensional BSDEs with jumps combining
the approach used in Hu-Peng [10].

This paper is organized as follows: In the next section, we state
some basic assumptions and basic estimates for BSDEs with jumps. And
then we state the main result of the paper on the BSVP. In Section
3, we apply our main result to the comparison theorem for BSDEs with
jumps. Finally, in Appendix, we give the proof of our main result.

\section{The BSVP with jumps in closed sets}\label{sec:intro}

\qquad Let $(\Omega,{\cal{F}},P,({\cal{F}}_{t})_{t\geq 0})$ be a
complete stochastic basis such that $\mathcal{F}_{0}$ contains all
$P$-null elements of ${\cal{F}}$, and
$\mathcal{F}_{t^{+}}:=\cap_{\ep>0}\mathcal{F}_{t+\ep}=\mathcal{F}_{t},t\geq
0$, and $\mathcal{F}=\mathcal{F}_{T}$, and suppose that the
filtration is generated by the following
two mutually independent processes:\\
(i) a $d$-dimensional standard Brownian motion $(W_{t})_{0\leq
t\leq T}$, and\\
(ii) a stationary Poisson random measure $N$ on $(0,T]\times E$,
where $E\subset R^{l}\setminus\{0\}$, $E$ is equipped with its Borel
field  $\mathcal {B}_E$, with compensator $\hat{N}(dtde)=dtn(de)$,
such that $n(E)<\infty$, and $\{\tilde{N}((0,t]\times
A)=(N-\hat{N})((0,t]\times A)\}_{0<t\leq T}$ is an
$\mathcal{F}_{t}$-martingale, for each $A\in \mathcal{B}_{E}$.

By $T>0$ we denote the finite real time horizon. We define some
spaces of processes. Let $\mathcal{S}^{2}_{[0,T]}$ denote the set of
$\mathcal{F}_{t}$-adapted c$\grave{a}$dl$\grave{a}$g $m$-dimensional
processes $\{Y_{t},0\leq t\leq T\}$ which are such that
$$
\|Y\|_{\mathcal{S}^{2}_{[0,T]}}:=(E[\sup_{0\leq t\leq
T}|Y_{t}|^{2}])^{\frac{1}{2}}<\infty.
$$
Let $L^{2}_{[0,T]}(W)$ be the set of $\mathcal{F}$-progressively
measurable $m\times d$ dimensional processes $\{Z_{t},0\leq t\leq
T\}$ which are such that
$$
\|Z\|_{L^{2}_{[0,T]}(W)}:=(E\displaystyle\int_{0}^{T}|Z_{t}|^{2}dt)^{\frac{1}{2}}<\infty.
$$
By $L^{2}_{[0,T]}(\tilde{N})$ we denote the set of mapping $U :
\Omega\times [0,T]\times E\rightarrow R^{m}$ which are
$\mathcal{P}\times \mathcal{B}_{E}$ measurable and such that
$$
\|U\|_{L^{2}_{[0,T]}(\tilde{N})}:=(E\displaystyle\int_{0}^{T}\displaystyle\int_{E}|U_{t}(e)|^{2}n(de)dt)^{\frac{1}{2}}<\infty,
$$
where $\mathcal{P}$ denotes the $\sigma$-algebra of
$\mathcal{F}_{t}$-predictable subsets of $\Omega\times [0,T]$.
Finally we define
$$\mathcal{B}^{2}_{[0,T]}=\mathcal{S}^{2}_{[0,T]}\times
L^{2}_{[0,T]}(W)\times L^{2}_{[0,T]}(\tilde{N}).$$

Consider the following BSDE with jump:
$$
Y_{t}=\xi+\displaystyle\int_{t}^{T}f(s,Y_{s},Z_{s},U_{s})ds-\displaystyle\int_{t}^{T}Z_{s}dW_{s}-\displaystyle\int_{t}^{T}\displaystyle\int_{E}U_{s}(e)\tilde{N}(dsde),0\leq
t\leq T, \eqno{(2.1)}
$$
where $\xi\in L^{2}(\Omega,{\cal{F}}_{T},P,R^{m})$ is a given random
variable and $f: \Omega\times [0,T]\times R^{m}\times R^{m\times
d}\times L^{2}(E,\mathcal{B}_{E},n; R^{m})\rightarrow R^{m}$ a
measurable function.

We suppose that there exists $L>0$ such that
$$
\begin{array}{ll}
(H1)& f(\cdot,\cdot,y,z,u) \mbox{ \ is progressively measurable, \
}f(\omega,\cdot,y,z,u) \mbox{\ is continuous},\\
(H2) & |f(t,y,z,u)-f(t,y',z',u')|\leq L(|y-y'|+|z-z'|+\|u-u'\|),\\
(H3)& \sup_{0\leq t\leq T}|f(t,0,0,0)|\in
L^{2}(\Omega,\mathcal{F}_{T},P,R^{m}),
\end{array}
$$
for all $0\leq t\leq T$, $y,y'\in R^{m}$, $z,z'\in R^{m\times d}$,
$u,u'\in L^{2}(E,\mathcal{B}_{E},n; R^{m})$, $P-a.s.$.

Let's recall the existence and uniqueness result for BSDEs with
jumps (see [5]):

{\bf Proposition 2.1.} {\it Let (H1)-(H3) holds true. Then for any
given\\ $\xi\in L^{2}(\Omega,{\cal{F}}_{T},P,R^{m})$, there exists
an unique triple $(Y,Z,U)\in \mathcal{B}^{2}_{[0,T]}$ which solves
BSDE (2.1).}

Note that the notion of BSDE with jump generalizes the well-known
martingale representation property. Indeed, in the particular case
of $f=0$ we have the following Lemma (see [11]):

{\bf Lemma 2.2. }{\it For any $\xi\in
L^{2}(\Omega,{\cal{F}}_{T},P,R^{m})$, there exist a unique $R(\xi)$
belonging to $L^{2}_{[0,T]}(W)$ and a $V\in
L^{2}_{[0,T]}(\tilde{N})$ such that
$$
\xi=E(\xi)+\displaystyle\int_{0}^{T}R_{s}dW_{s}+\displaystyle\int_{0}^{T}\displaystyle\int_{E}V_{s}(e)\tilde{N}(dsde).
\eqno{(2.2)}
$$
}

 For completeness, we give the following definition of stochastic
viability:

{\bf Definition 2.3.} {\it Let $K$ be a nonempty closed subset of $R^{m}$.\\
{\rm(a)} A stochastic process $\{Y_{t},t\in [0,T]\}$ is viable in
$K$ if and only if for $P$-almost $\omega\in \Omega$,
$$
Y_{t}(\omega)\in K,\mbox{ \ \ }\forall t\in [0,T].
$$
{\rm(b)} The closed subset $K$ enjoys the BSVP for the equation
(2.1) if and only if:

$\forall t\in [0,T]$, $\forall \xi\in
L^{2}(\Omega,{\cal{F}}_{t},P,K)$, there exists a solution
$(Y,Z,U)\in \mathcal{B}^{2}_{[0,t]}$ to BSDE (2.1) over the time
interval $[0,t]$,
$$
Y_{s}=\xi+\displaystyle\int_{s}^{t}f(r,Y_{r},Z_{r},U_{r})dr-\displaystyle\int_{t}^{T}Z_{r}dW_{r}-\displaystyle\int_{t}^{T}\displaystyle\int_{E}U_{r}(e)\tilde{N}(drde),s\in
[0,t],
$$
such that $\{Y_{s},s\in [0,t]\}$ is viable in $K$.}

The question we are interested in is that when $K$ enjoys the BSVP
for BSDE (2.1).

Let us define for any closed set $K\subset R^{m}$ the multivalued
projection of a point $a$ onto $K$:
$$
\Pi_{K}(a):=\{b\in K | \mbox{ \ }|a-b|=\min_{c\in
K}|a-c|=d_{K}(a)\}.
$$

Recall that $\Pi_{K}(a)$ is a singleton whenever $d_{K}$ is
differentiable at the point $a$. According to Motzkin's Theorem
$\Pi_{K}$ is single-valued if and only if $K$ is convex. Notice that
$d^{2}_{K}(\cdot)$ is convex when $K$ is convex, and thus, due to
Alexandrov's Theorem [1], $d^{2}_{K}(\cdot)$ is almost everywhere
twice differentiable. By twice differentiable, we mean that the
function admits a second order Taylor expansion. And this may hold
true even if the first derivative is not continuous.

We need an auxiliary result on BSDEs with jumps. Give any $\xi\in
L^{2}(\Omega,{\cal{F}}_{T},P,R^{m})$, we denote by $(Y,Z,U)$ the
unique solution to BSDE (2.1) and by $R$, $V$ the processes
associated to $\xi$ by Lemma 2.2. With these notations we can state

{\bf Proposition 2.4.} {\it Suppose that (H1)-(H3) hold true. Then
there exist real positive constants $C_{0},M$ such that for all
$\xi\in L^{2}(\Omega,{\cal{F}}_{T},P,R^{m})$,
$$
\begin{array}{ll}
& E[\sup_{s\in
[t,T]}|Y_{s}|^{2}|\mathcal{F}_{t}]+E[\displaystyle\displaystyle\int_{t}^{T}|Z_{s}|^{2}ds|\mathcal{F}_{t}]+E[\displaystyle\displaystyle\int_{t}^{T}\displaystyle\displaystyle\int_{E}|U_{s}(e)|^{2}n(de)ds|\mathcal{F}_{t}]\\
\leq &
C_{0}\{E[|\xi|^{2}|\mathcal{F}_{t}]+E[(\displaystyle\displaystyle\int_{t}^{T}|f(s,0,0,0)|ds)^{2}|\mathcal{F}_{t}]\},t\in
[0,T],
\end{array}
\eqno{(2.3)}
$$
and for all $\ep\in [0,T]$,
$$
\begin{array}{ll}
& E[\sup_{s\in
[T-\ep,T]}|Y_{s}-E(\xi|\mathcal{F}_{s})|^{2}|\mathcal{F}_{T-\ep}]
+E[\displaystyle\displaystyle\int_{T-\ep}^{T}|Z_{s}-R_{s}|^{2}ds|\mathcal{F}_{T-\ep}]\\
&+E[\displaystyle\displaystyle\int_{T-\ep}^{T}\displaystyle\displaystyle\int_{E}|U_{s}(e)-V_{s}(e)|^{2}n(de)ds|\mathcal{F}_{T-\ep}]\\
\leq & C_{0}\ep
E[\displaystyle\displaystyle\int_{T-\ep}^{T}|f(s,E(\xi|\mathcal{F}_{s}),R_{s},V_{s})|^{2}ds|\mathcal{F}_{T-\ep}],
\end{array}
\eqno{(2.4)}
$$
in particular, there exists a positive constant $M$ such that
$$
\begin{array}{ll}
&E[|Y_{t}-E(\xi|\mathcal{F}_{t})|^{2}]+E[\displaystyle\int_{t}^{T}|Z_{s}-R_{s}|^{2}ds]+E[\displaystyle\int_{t}^{T}\displaystyle\int_{E}|U_{s}(e)-V_{s}(e)|^{2}n(de)ds]\\
\leq& (T-t)M,t\in [0,T].\end{array}\eqno{(2.5)}
$$}

Now we can get the following theorem.

{\bf Theorem 2.5.}  {\it Suppose that $f : \Omega\times [0,T]\times
R^{m}\times R^{m\times d}\times L^{2}(E,\mathcal{B}_{E},n;
R^{m})\rightarrow R^{m}$ is a measurable function which satisfies
(H1)-(H3). Let $K$ be a nonempty closed set. If $K$ enjoys the BSVP
for BSDE (2.1), then the set $K$ is convex.}

{\bf Proof: }The method is totally the same as that in Theorem 2.4
in [7] thanks to (2.5) in Proposition 2.4. So we omit it.

Since the previous theorem means that only convex sets could have
the BSVP even with jumps, we restrict our attention to closed convex
sets.

{\bf Theorem 2.6.} {\it Suppose that $f : \Omega\times [0,T]\times
R^{m}\times R^{m\times d}\times L^{2}(E,\mathcal{B}_{E},n;
R^{m})\rightarrow R^{m}$ is a measurable function which satisfies
(H1)-(H3). Let $K$ be a nonempty closed set. The set $K$ enjoys the
BSVP for BSDE (2.1) if and only if:

$\forall (t,z,u)\in [0,T]\times R^{m\times d}\times
L^{2}(E,\mathcal{B}_{E},n; R^{m})$ and for all $y\in R^{m}$ such
that $d^{2}_{K}(\cdot)$ is twice differentiable at $y$,
$$
\begin{array}{ll}
&4\langle y-\Pi_{K}(y),f(t,y,z,u)\rangle \\
\leq &\langle D^{2}d^{2}_{K}(y)z,z\rangle+
C^{*}d^{2}_{K}(y)\\
&+2\displaystyle\int_{E}[d^{2}_{K}(y+u(e))-d^{2}_{K}(y)-2\langle
y-\Pi_{K}(y),u(e)\rangle ]n(de),P-a.s. ,
\end{array}\eqno(2.6)
$$
where $C^{*}>0$ is a constant which does not depend on $(t,y,z,u)$.}

Let us notice that, under assumption (H2), condition (2.6) takes
form:
$$
\begin{array}{ll}
&4\langle y-\Pi_{K}(y),f(t,\Pi_{K}(y),z,u)\rangle \\
\leq &\langle D^{2}d^{2}_{K}(y)z,z\rangle+(
C^{*}+4L)d^{2}_{K}(y)\\
&+2\displaystyle\int_{E}[d^{2}_{K}(y+u(e))-d^{2}_{K}(y)-2\langle
y-\Pi_{K}(y),u(e)\rangle ]n(de),P-a.s..
\end{array}
$$
On the other hand, for some $C'>0$, the condition
$$
\begin{array}{ll}
&4\langle y-\Pi_{K}(y),f(t,\Pi_{K}(y),z,u)\rangle \\
\leq &\langle
D^{2}d^{2}_{K}(y)z,z\rangle+C'd^{2}_{K}(y)+\\
&2\displaystyle\int_{E}[d^{2}_{K}(y+u(e))-d^{2}_{K}(y)-2\langle
y-\Pi_{K}(y),u(e)\rangle ]n(de),P-a.s.
\end{array}\eqno{(2.7)}
$$
implies (2.6) with constant $C'+4L$ instead of $C^{*}$.

This shows that (2.6) is a condition only on the values of
$f(t,\cdot,z,u)$ on $\partial K$. Recall also that the behavior of
$d_{K}$ on $R^{m}$ is completely determined by that of $\partial K$.

{\bf Remark 2.7.} {\it Because $D^{2}[d^{2}_{K}]$ is almost
everywhere positive semidefinite,  and
$$
\begin{array}{ll}
&d^{2}_{K}(y+u(e))-d^{2}_{K}(y)-2\langle
y-\Pi_{K}(y),u(e)\rangle\\
= &
(y+u(e)-\Pi_{K}(y)+\Pi_{K}(y)-\Pi_{K}(y+u(e)))^{2}-(y-\Pi_{K}(y))^{2}-2\langle
y-\Pi_{K}(y),u(e)\rangle\\
=& |u(e)+\Pi_{K}(y)-\Pi_{K}(y+u(e))|^{2}+2\langle
y-\Pi_{K}(y),\Pi_{K}(y)-\Pi_{K}(y+u(e))\rangle\\
\geq & 0.
\end{array}
$$
The last inequality is from the convexity of $K$. So a sufficient
condition for the BSVP of $K$ with jumps is that there exists some
constant $C^{*}>0$, such that, for all $(t,y,z,u)$,
$$
4\langle y-\Pi_{K}(y),f(t,y,z,u)\rangle\leq
C^{*}d^{2}_{K}(y),P-a.s.. \eqno{(2.7')}
$$
}

{\bf Example 2.8.} {\it If we set
$K=\{(x_{1},x_{2})|x^{2}_{1}+x^{2}_{2}\leq 1\}$. In this case,
$$
d^{2}_{K}(x)=\left\{
\begin{array}{ll}
0,&\mbox{when \ }x\in K,\\
(|x|-1)^{2},& \mbox{when \ }x\in R^{2}\backslash K.
\end{array}\right.
$$
Moreover, $d^{2}_{K}(\cdot)\in C^{1}(R^{2})$ and twice
differentiable when $x\in R^{2}\backslash \partial K$. To satisfy
$(2.7')$, we can choose $$f(t,y,z,u)=y-\Pi_{K}(y).$$ So obviously,
$(2.7')$ holds true with $C^{*}\geq 4$. In fact the following BSDE
with jump does enjoy BSVP for $K$:
$$
Y_{t}=\xi+\displaystyle\int_{t}^{T}(Y_{s}-\Pi_{K}(Y_{s}))ds-\displaystyle\int_{t}^{T}Z_{s}dW_{s}-\displaystyle\int_{t}^{T}\displaystyle\int_{E}U_{s}(e)\tilde{N}(dsde),0\leq
t\leq T.
$$
In deed, the following BSDE with jump has the same solution with
above when $\xi\in K$:
$$
Y_{t}=\xi-\displaystyle\int_{t}^{T}Z_{s}dW_{s}-\displaystyle\int_{t}^{T}\displaystyle\int_{E}U_{s}(e)\tilde{N}(dsde),0\leq
t\leq T.
$$
So $$|Y_{t}|^{2}=|E[\xi|\mathcal{F}_{t}]|^{2}\leq
E[|\xi|^{2}|\mathcal{F}_{t}]\leq 1.$$}

\section{The comparison theorem of multidimensional BSDEs with jumps}\label{sec:intro}
\qquad In this section, with the BSVP theory with jumps we have got
in the previous section, we can use the same method as in [10] to
study the comparison theorem for multidimensional BSDEs with jumps.

Consider the following BSDEs: $i=1,2,$
$$
Y^{i}_{t}=\xi^{i}+\displaystyle\int_{t}^{T}f^{i}(s,Y^{i}_{s},Z^{i}_{s},U^{i}_{s})ds-\displaystyle\int_{t}^{T}Z^{i}_{s}dW_{s}-\displaystyle\int_{t}^{T}\displaystyle\int_{E}U^{i}_{s}\tilde{N}(dsde),
\eqno{(3.1)}
$$
where $f^{1},f^{2}$ satisfy (H1)-(H3), and $\xi^{1},\xi^{2}\in
L^{2}(\Omega,\mathcal{F},P)$. In this subsection, we study the
following problem: under which condition the comparison theorem
holds? Interestingly, this problem is transformed to a viability
problem in $R^{m}_{+}\times R^{m}$ of $(Y^{1}-Y^{2},Y^{2})$.

{\bf Theorem 3.1.} {\it Suppose that $f^{1}$ and $f^{2}$ satisfy
(H1)-(H3). Then the following are equivalent:

{\rm(i)} For any $s\in [0,T]$, $\forall \xi^{1},\xi^{2}\in
L^{2}(\Omega,\mathcal{F}_{s},P;R^{m})$ such that $\xi^{1}\geq
\xi^{2}$, the unique solutions $(Y^{1},Z^{1},U^{1})$ and
$(Y^{2},Z^{2},U^{2})$ in $\mathcal{B}^{2}_{[0,s]}$ to the BSDE (3.1)
over interval $[0,s]$ satisfy:
$$
Y^{1}_{t}\geq Y^{2}_{t},\mbox{ \ }t\in [0,s];
$$

{\rm(ii)}$\forall t\in [0,T]$, $\forall (y,y')$, $\forall (z,z')$,
$\forall (u,u')$,
$$
\begin{array}{ll}
&-4\langle y^{-},f^{1}(t,y^{+}+y',z,u)-f^{2}(t,y',z',u')\rangle\\
\leq &
2\sum\limits_{k=1}^{m}I_{\{y_{k}<0\}}|z_{k}-z'_{k}|^{2}+C|y^{-}|^{2}+2\sum\limits_{k=1}^{m}I_{\{
y_{k}\geq
0\}}\displaystyle\int_{E}|(y_{k}+u_{k}(e)-u'_{k}(e))^{-}|^{2}n(de)\\
&+2\sum\limits_{k=1}^{m}I_{\{y_{k}<0\}}\displaystyle\int_{E}[|(y_{k}+u_{k}(e)-u'_{k}(e))^{-}|^{2}-|y_{k}^{-}|^{2}
-2y_{k}(u_{k}(e)-u'_{k}(e))]n(de),P-a.s.,\\
 \end{array}\eqno{(3.2)}
$$
where $C$ is a constant which dose not depend on
$t,(y,y'),(z,z'),(u,u')$. }

{\bf Proof: }Set
$$
\bar{Y}_{t}=(Y^{1}_{t}-Y^{2}_{t},Y^{2}_{t}),\bar{Z}_{t}=(Z^{1}_{t}-Z^{2}_{t},Z^{2}_{t}),\bar{U}_{t}=(U^{1}_{t}-U^{2}_{t},U^{2}_{t}).
$$
Then (i) is equivalent to the following:

(iii) For any $s\in [0,T]$,$\forall
\bar{\xi}=(\bar{\xi}^{1},\bar{\xi}^{2})$ such that
$\bar{\xi}^{1}\geq 0$, the unique solution
$(\bar{Y},\bar{Z},\bar{U})$ to the following BSDE over time interval
$[0,s]$:
$$
\bar{Y}_{t}=\bar{\xi}+\displaystyle\int_{t}^{s}\bar{f}(s,\bar{Y}_{s},\bar{Z}_{s},\bar{U}_{s})ds-\displaystyle\int_{t}^{s}\bar{Z}_{s}dW_{s}-\displaystyle\int_{t}^{s}\displaystyle\int_{E}\bar{U}_{s}\tilde{N}(dsde),
\eqno{(3.3)}
$$
satisfies $\bar{Y}^{1}\geq 0$, where for
$\bar{y}=(\bar{y}^{1},\bar{y}^{2})$,
$\bar{z}=(\bar{z}^{1},\bar{z}^{2})$,
$\bar{u}=(\bar{u}^{1},\bar{u}^{2})$,
$$
\bar{f}(t,\bar{y},\bar{z},\bar{u})=(f^{1}(t,,\bar{y}^{1}+\bar{y}^{2},\bar{z}^{1}+\bar{z}^{2},\bar{u}^{1}+\bar{u}^{2})
-f^{2}(t,\bar{y}^{2},\bar{z}^{2},\bar{u}^{2}),f^{2}(t,\bar{y}^{2},\bar{z}^{2},\bar{u}^{2})).
$$

So we can apply Theorem 2.6 to BSDE (3.3) and the convex closed set
$K:=R^{m}_{+}\times R^{m}$. Obviously, when $\hat{y}=(y,y')\in K$,
(2.7) and (3.2) holds true naturally. When $\hat{y}=(y,y')\in
R^{2m}\diagdown K$, since $\forall x=(x^1,x^2)\in R^{2m}$,
$$
\Pi_{K}(x)=\left(\begin{array}{c}
(x^{1})^{+}\\\\
x^{2}
\end{array}\right),
x-\Pi_{K}(x)=\left(\begin{array}{c}
-(x^{1})^{-}\\\\
0
\end{array}\right).
$$
So
$$
d^{2}_{K}(x)=|(x^{1})^{-}|^{2}=\sum\limits_{k=1}^{m}I_{\{x^{1}_{k}<0\}}|x^{1}_{k}|^{2}.
$$
And
$$
(D^{2}d^{2}_{K})(x)\left\{\begin{array}{ll} =0_{2m\times 2m},&
\mbox{when \ }x\in K^\circ,\\
\mbox{does not exist},& \mbox{when \ }x\in \partial K,\\
=(a_{ij})_{2m\times 2m},& \mbox{when \ }x\in R^{2m}\backslash K,
\end{array}
\right.
$$
where
$$
a_{ij}=0,\mbox{ \ when \ }i\neq j,\mbox{ \ } a_{ii}=\left\{
\begin{array}{ll}
0,& m<i\leq 2m,\\
0,& 1\leq i\leq m,x^{1}_{i}\geq 0,\\
2,& 1\leq i\leq m,x^{1}_{i}< 0.
\end{array}
\right.
$$
We can get that (2.7) is equivalent to
$$
\begin{array}{ll}
& -4\langle y^{-},f^{1}(t,y^{+}+y',z,u)-f^{2}(t,y',z',u')\rangle\\
\leq &
2\sum\limits_{k=1}^{m}I_{\{y_{k}<0\}}|z_{k}-z'_{k}|^{2}+C|y^{-}|^{2}
+2\sum\limits_{k=1}^{m}I_{y_{k}\geq 0}\displaystyle\int_{E}|(y_{k}+u_{k}(e)-u'_{k}(e))^{-}|^{2}n(de)\\
&+2\sum\limits_{k=1}^{m}I_{\{y_{k}<0\}}\displaystyle\int_{E}[|(y_{k}+u_{k}(e)-u'_{k}(e))^{-}|^{2}-|y_{k}|^{2}-2
y_{k}(u_{k}(e)-u'_{k}(e))]n(de).
\end{array}
$$
$$
\eqno{\Box}
$$

{\bf Theorem 3.2.} {\it Let $m=1$. Suppose that $f^{1}$ and $f^{2}$
satisfy (H1)-(H3). Then the following are equivalent:

{\rm(i)} For any $s\in [0,T]$, $\forall \xi^{1},\xi^{2}\in
L^{2}(\Omega,\mathcal{F}_{s},P;R)$ such that $\xi^{1}\geq \xi^{2}$,
the unique solutions $(Y^{1},Z^{1},U^{1})$ and $(Y^{2},Z^{2},U^{2})$
in $\mathcal{B}^{2}_{[0,s]}$ to the BSDE (3.1) over interval $[0,s]$
satisfy:
$$
Y^{1}_{t}\geq Y^{2}_{t},\mbox{ \ }t\in [0,s];
$$

{\rm(ii)} $\forall (t,y',z)\in [0,T]\times R\times R^{1\times d}$,
$$
\begin{array}{l}
\forall u,u'\in L^{2}(E,\mathcal{B}_{E},n;R) \mbox{ \ such
that \ }u\geq u', n(de)-a.s.,\\
f^{1}(t,y',z,u)-f^{2}(t,y',z,u')\geq
-\displaystyle\int_{E}(u(e)-u'(e))n(de).
\end{array}
 \eqno{(3.4)}
$$
} {\bf Proof: }(i)$\Rightarrow$ (ii). Since $m=1$, choose
$y<0,z=z',u\geq u'$, from (3.2), we have
$$
\begin{array}{ll}
&4y(f^{1}(t,y',z',u)-f^{2}(t,y',z',u'))\\
\leq&
Cy^{2}+2\displaystyle\int_{E}[|(y+u(e)-u'(e))^{-}|^{2}-y^{2}-2y(u(e)-u'(e))]n(de)\\
\leq & Cy^{2}+2\displaystyle\int_{E}-2y(u(e)-u'(e))n(de).
\end{array}
$$
Dividing by $4y$ and letting $y$ tend to 0, we get (3.4).

(ii)$\Rightarrow$ (i). When $y\geq 0$, (3.2) holds true naturally.
When $y<0$, from (3.4) we have: $\forall t\in [0,T],y',z,z',u,u',$
$$
\begin{array}{ll}
& 4y [f^{1}(t,y',z,u)-f^{2}(t,y',z',u')]\\
=&
4y[f^{1}(t,y',z,u)-f^{1}(t,y',z',u)]+4y[f^{1}(t,y',z',u)-f^{2}(t,y',z',I_{\{u>
u'\}}u'+I_{\{u\leq u'\}}u)]\\
&+4y[f^{2}(t,y',z',I_{\{u>
u'\}}u'+I_{\{u\leq u'\}}u)-f^{2}(t,y',z',u')]\\
\leq&
2L^{2}y^{2}+\frac{2}{L^{2}}|f^{1}(t,y',z,u)-f^{1}(t,y',z',u)|^{2}+\displaystyle\int_{\{u>u'\}}-4y(u(e)-u'(e))n(de)\\
&+2L^{2}y^{2}+\frac{2}{L^{2}}|f^{2}(t,y',z',I_{\{u>
u'\}}u'+I_{\{u\leq u'\}}u)-f^{2}(t,y',z',u')|^{2}\\
\leq &(4L^{2}+2n(E))y^{2}+2|z-z'|^{2} +2\displaystyle\int_{\{u\leq
u'\}}|u(e)-u'(e)|^{2}n(de)\\
& +\displaystyle\int_{\{u>u'\}}(-2y^{2})n(de)+\displaystyle\int_{\{u>u'\}}-4y(u(e)-u'(e))n(de)\\
\leq &(4L^{2}+2n(E))y^{2}+2|z-z'|^{2} +2\displaystyle\int_{\{u\leq
u'\}}|u(e)-u'(e)|^{2}n(de)\\
&
+2\displaystyle\int_{\{u>u'\}}[|(y+u(e)-u'(e))^{-}|^{2}-y^{2}-2y(u(e)-u'(e))]n(de)\\
=& (4L^{2}+n(E))y^{2}+2|z-z'|^{2}
+2\displaystyle\int_{E}[|(y+u(e)-u'(e))^{-}|^{2}-y^{2}-2y(u(e)-u'(e))]n(de)
 \end{array}
$$
$$
\eqno{\Box}
$$

{\bf Remark 3.3.} {\it When $m=1$, (3.4) constrains the form of the
dependence of the generators in the jump components of the BSDEs. We
refer the reader to Wu [15], Royer [14] for the analogue of (3.4)
which can be somewhat contained by (3.4). But the most important of
all, (3.4) is a necessary and sufficient condition for the
comparison theorem to hold. While in almost any other references on
the comparison theorem for 1-dimensional BSDEs with jumps, the
results are just sufficient. }

{\bf Remark 3.4.} {\it The following example and counter-example
show that (3.4) is really a necessary and sufficient condition for
comparison theorem.

Let $E=R\setminus \{0\}$, $n(de)=\d_1(de)$. Then
$$N_t=\displaystyle\int_0^t\displaystyle\int_E N(dsde), 0\leq t\leq T,$$
is a standard Poisson process.

(a) If we choose $$f^1(t,y,z,u)=f^2(t,y,z,u)=-\frac{1}{2}u(1),$$
then $\forall (t,y,z)\in [0,T]\times R\times R^{1\times d}$ and
$\forall u,u'\in L^{2}(E,\mathcal{B}_{E},n;R)$ such that $$u\geq u',
 n(de)-a.s.,$$ we have
$$
f(t,y,z,u)-f(t,y,z,u')=-\frac{1}{2}(u(1)-u'(1))\geq
-(u(1)-u'(1))=-\displaystyle\int_E (u(e)-u'(e))n(de).
$$
This means that (3.4) holds true. So We know that comparison theorem
holds. For example, for any $s\in [0,T]$, if we choose
$$
\xi^1=N_s\geq \xi^2=0,
$$
then for all $t\in [0,s],$
$$
\begin{array}{l}
(Y^{1}_t,Z^1_t,U^1_t)=(N_t+\frac{1}{2}(s-t),0,I_{\{e=1\}}),\\
(Y^{2}_t,Z^2_t,U^2_t)=(0,0,0).
\end{array}
$$
It's clear that $P\{Y^{1}_t\geq Y^{2}_t\}=1$, for all $0\leq t\leq
s$.

(b)On the other hand, if we choose
$$f^1(t,y,z,u)=f^2(t,y,z,u)=-2u(1)$$ and
$$
\xi^1=N_s, \xi^2=0,\forall s\in [0,T],
$$
then for all $t\in [0,s],$
$$
\begin{array}{l}
(Y^{1}_t,Z^1_t,U^1_t)=(N_t-(s-t),0,I_{\{e=1\}}),\\
(Y^{2}_t,Z^2_t,U^2_t)=(0,0,0).
\end{array}
$$
Obviously, $N_s\geq 0, P-a.s.$, but $P\{Y^{1}_t< Y^{2}_t\}>0$, for
all $0\leq t\leq s$. This is because, in this case,
$$
f(t,y,z,u)-f(t,y,z,u')=-2(u(1)-u'(1))<
-(u(1)-u'(1))=-\displaystyle\int_E (u(e)-u'(e))n(de).
$$
}

{\bf Corollary 3.5.} {\it Let $m=1$ and  suppose furthermore that
$f^{1}$ and $f^{2}$ are independent of $u$ and satisfy (H1)-(H3).
Then the following are equivalent:

{\rm(i)} For any $s\in [0,T]$, $\forall \xi^{1},\xi^{2}\in
L^{2}(\Omega,\mathcal{F}_{s},P;R)$ such that $\xi^{1}\geq \xi^{2}$,
the unique solutions $(Y^{1},Z^{1},U^{1})$ and $(Y^{2},Z^{2},U^{2})$
in $\mathcal{B}^{2}_{[0,s]}$ to the BSDE (3.1) over interval $[0,s]$
satisfy:
$$
Y^{1}_{t}\geq Y^{2}_{t},\mbox{ \ }t\in [0,s];
$$

{\rm(ii)}$ f^{1}(t,y,z)\geq f^{2}(t,y,z),\forall t\in [0,T],\forall
(y,z)\in R\times R^{1\times d}. $}

 This generalizes the result in [10].

Now consider the special case when $f^{1}=f^{2}=f$. Then we have:

{\bf Theorem 3.6.} {\it Suppose that $f$ satisfies (H1)-(H3). Then
the following are equivalent:

{\rm(i)} For any $s\in [0,T]$, $\forall \xi^{1},\xi^{2}\in
L^{2}(\Omega,\mathcal{F}_{s},P;R^{m})$ such that $\xi^{1}\geq
\xi^{2}$, the unique solutions $(Y^{1},Z^{1},U^{1})$ and
$(Y^{2},Z^{2},U^{2})$ in $\mathcal{B}^{2}_{[0,s]}$ to the BSDE (3.1)
over interval $[0,s]$ satisfy:
$$
Y^{1}_{t}\geq Y^{2}_{t},\mbox{ \ }t\in [0,s];
$$

{\rm(ii)} For any $k=1,2,...,m,$
$$
\left\{
\begin{array}{ll}
{\rm(a)}& f_{k} \mbox{ \ depends only on \ } z_{k};\\
{\rm(b)}& f_{k}(t,\d^{k}y+y',z_{k},u)-f_{k}(t,y',z_{k},u')\geq
-\displaystyle\int_{E}(u_{k}(e)-u'_{k}(e))n(de),\\
& \forall u\geq u',n(de)-a.s., \mbox{ \ for any \ }\d^{k}y\in R^{m}
\mbox{ \ such that \ }\d^{k}y\geq 0,(\d^{k}y)_{k}=0;\\
{\rm(c)}& \mbox{There exists a positive constant \ } C_{0},\mbox{ \ s.t., \ }\\
&\forall y,y',z,u,u' \mbox{ \ satisfy \ }u\leq
u',n(de)-a.s.,\\
& \sum\limits_{k=1}^{m}I_{\{y_{k<0}\}}4y_{k}
(f_{k}(t,y^{+}+y',z_{k},u)-f_{k}(t,y',z_{k},u'))\\
&\leq
C_{0}|y^-|^2+2\sum\limits_{k=1}^{m}I_{\{y_{k<0}\}}\displaystyle\int_{E}|u_{k}(e)-u'_{k}(e)|^{2}n(de)\\
&\mbox{ \
}+2\sum\limits_{k=1}^{m}I_{\{y_{k\geq0}\}}\displaystyle\int_{E}|(y_{k}+u_{k}(e)-u'_{k}(e))^{-}|^{2}n(de).
\end{array}
\right.\eqno{(3.5)}
$$
}

 {\bf Proof: }From Theorem 3.1, (i) is equivalent to (3.2) with
$f^{1}=f^{2}$. What we have to do is to prove that, when
$f^{1}=f^{2}$, (3.2) $\Leftrightarrow$ (ii).

Let us pick $y_{k}<0$, and $y=y_{k}e_{k}$, $z_{k}=z'_{k}$, $u=u'$,
from (3.2), we have
$$
4y_{k}(f_{k}(t,y',z,u)-f_{k}(t,y',z',u))\leq Cy^{2}_{k}.
$$
We deduce easily that
$$
|f_{k}(t,y',z,u)-f_{k}(t,y',z',u)|=0, \mbox{ \ when \ }z_{k}=z'_{k}.
$$
So obviously, $f_{k}$ depends only on $z_{k}$.

Moreover, for $\d^{k}y\in R^{m}$ such that $\d^{k}y\geq
0,(\d^{k}y)_{k}=0$, putting in (3.2) with $y=\d^{k}y-\ep
e_{k},\ep>0$, $z=z'$, $u(e)\geq u'(e), n(de)-a.s.$,  we have
$$
\begin{array}{ll}
&-4\ep(f_{k}(t,\d^{k}y+y',z_{k},u)-f_{k}(t,y',z_{k},u'))\\
\leq&
C\ep^{2}+2\displaystyle\int_{E}[|(-\ep+u_{k}(e)-u'_{k}(e))^-|^{2}-\ep^{2}+2\ep(u_{k}(e)-u'_{k}(e))]n(de).
\end{array}
$$
Dividing by $-4\ep$ and letting $\ep\rightarrow 0$, we hence get (b)
in (3.5).

For (c), it's straightforward by (3.2) with $u\leq u',n(de)-a.s..$

(ii)$\Rightarrow$(3.2). When $y\geq 0$, (3.2) holds true naturally.
When there exist $1\leq k\leq n$, such that, $y_{k}<0$. $\forall
t\in [0,T],y',z,z',u,u'$, we set $\hat{u}=(\hat{u}_{k})_{k=1}^{m}$,
where
$$\hat{u}_{k}=I_{\{u_{k}>u'_{k}\}}u'_{k}+I_{\{u_{k}\leq
u'_{k}\}}u_{k}.$$ So we have $u\geq \hat{u}$, $\hat{u}\leq u'$. Then
from (ii) we have:
$$
\begin{array}{ll}
&\sum\limits_{k=1}^{m}I_{\{y_{k}<0\}}4y_{k}[f_{k}(t,y^{+}+y',z_{k},u)-f_{k}(t,y',z'_{k},u')]\\
=&\sum\limits_{k=1}^{m}I_{\{y_{k}<0\}}4y_{k}[f_{k}(t,y^{+}+y',z_{k},u)-f_{k}(t,y^{+}+y',z'_{k},u)\\
&\mbox{\ \ \ \ \ \ \ \ \ \ \ \ \ \
 \ \ \ \  \ \ }+f_{k}(t,y^{+}+y',z'_{k},u)-f_{k}(t,y^{+}+y',z'_{k},\hat{u})\\
&\mbox{\ \ \ \ \ \ \ \ \ \ \ \ \ \
 \ \ \ \  \ \ }+f_{k}(t,y^{+}+y',z'_{k},\hat{u})-f_{k}(t,y',z'_{k},u')]\\
\leq
&\sum\limits_{k=1}^{m}I_{\{y_{k}<0\}}[(2L^{2}+2n(E)+C_{0})y^{2}_{k}+2|z_{k}-z'_{k}|^{2}+2\displaystyle\int_{\{u_{k}\leq
u'_{k}\}}|u_{k}(e)-u'_{k}(e)|^{2}n(de)]\\
&+\sum\limits_{k=1}^{m}I_{\{y_{k}<0\}}\displaystyle\int_{\{u_{k}>u'_{k}\}}2[|(y_{k}+u_{k}(e)-u'_{k}(e))^{-}|^{2}-y^{2}_{k}-2y_{k}(u_{k}(e)-u'_{k}(e))n(de)\\
&+2\sum\limits_{k=1}^{m}I_{\{y_{k}\geq
0\}}\displaystyle\int_{E}|(y_{k}+u_{k}(e)-u'_{k}(e))^{-}|^{2}n(de).
\end{array}
$$

{\bf Remark 3.7.} {\it When $f_{k}$ depends only on $u_{k}$, with
(b) in (3.5) and the Lipschitz condition of $f$ w.r.t. $u$ we have
$$\begin{array}{ll}
& \sum\limits_{k=1}^{m}I_{\{y_{k<0}\}}4y_{k}
(f_{k}(t,y^{+}+y',z_{k},u_k)-f_{k}(t,y',z_{k},u'_k))\\
=& \sum\limits_{k=1}^{m}I_{\{y_{k<0}\}}4y_{k}
(f_{k}(t,y^{+}+y',z_{k},u_k)-f_{k}(t,y',z_{k},u_k)+f_{k}(t,y',z_{k},u_k)-f_{k}(t,y',z_{k},u'_k))\\
\leq &\sum\limits_{k=1}^{m}I_{\{y_{k<0}\}}4y_{k}(
f_{k}(t,y',z_{k},u_k)-f_{k}(t,y',z_{k},u'_k))\\
\leq &
2L^2|y^-|^2+2\sum\limits_{k=1}^{m}I_{\{y_{k<0}\}}\displaystyle\int_E|u_k(e)-u'_k(e)|^2n(de).
\end{array}
$$
So the condition (c) in (3.5) can be cancelled:}

{\bf Theorem 3.8.} {\it Suppose that $f$ satisfies (H1)-(H3) and
$f_k$ depends only on $u_{k}$. Then the following are equivalent:

{\rm(i)} For any $s\in [0,T]$, $\forall \xi^{1},\xi^{2}\in
L^{2}(\Omega,\mathcal{F}_{s},P;R^{m})$ such that $\xi^{1}\geq
\xi^{2}$, the unique solutions $(Y^{1},Z^{1},U^{1})$ and
$(Y^{2},Z^{2},U^{2})$ in $\mathcal{B}^{2}_{[0,s]}$ to the BSDE (3.1)
over interval $[0,s]$ satisfy:
$$
Y^{1}_{t}\geq Y^{2}_{t},\mbox{ \ }t\in [0,s];
$$

{\rm(ii)} For any $k=1,2,...,m,$
$$
\left\{
\begin{array}{ll}
{\rm(a)}& f_{k} \mbox{ \ depends only on \ } z_{k};\\
{\rm(b)}& \forall u\geq u',n(de)-a.s., \mbox{ \ for any \
}\d^{k}y\in R^{m} \mbox{ \ such that \ }\d^{k}y\geq
0,(\d^{k}y)_{k}=0,\\
&f_{k}(t,\d^{k}y+y',z_{k},u_k)-f_{k}(t,y',z_{k},u'_k)\geq
-\displaystyle\int_{E}(u_{k}(e)-u'_{k}(e))n(de).
\end{array}\right.
$$
}

\vspace{3mm} {\bf Corollary 3.9.} {\it If we suppose that $f$ is
independent of $u$, then the following are equivalent:

{\rm(i)} For any $s\in [0,T]$, $\forall \xi^{1},\xi^{2}\in
L^{2}(\Omega,\mathcal{F}_{s},P;R^{m})$ such that $\xi^{1}\geq
\xi^{2}$, the unique solutions $(Y^{1},Z^{1},U^{1})$ and
$(Y^{2},Z^{2},U^{2})$ in $\mathcal{B}^{2}_{[0,s]}$ to the BSDE (3.1)
over interval $[0,s]$ satisfy:
$$
Y^{1}_{t}\geq Y^{2}_{t},\mbox{ \ }t\in [0,s];
$$

{\rm(ii)} For any $k=1,2,...,n,$
$$
f_{k} \mbox{ \ depends only on \ } z_{k} \mbox{ \ and}
$$
$$
f_{k}(t,\d^{k}y+y',z_{k})\geq f_{k}(t,y',z_{k}),\mbox{ \ for any \
}\d^{k}y\in R^{m} \mbox{ \ such that \ }\d^{k}y\geq
0,(\d^{k}y)_{k}=0.
$$ }
This generalizes the result in [10].

\vspace{3mm}

At the end of this section, we study the comparison theorem for
multidimensional BSDEs with jumps.

We define $\mathbb{S}^{m}$ as the space of symmetric real $m\times
m$ matrices, and denote by $\mathbb{S}_{+}^{m}$ the subspace of
$\mathbb{S}^{m}$ containing the nonnegative elements in
$\mathbb{S}^{m}$. Without loss of generality, we set $d=1$.

Let us again consider a function $F$, which will be in the sequel
the generator of the BSDE, defined on $\Omega\times [0,T]\times
\mathbb{S}^{m}\times \mathbb{S}^{m}\times \mathbb{S}^{m}$, with
values in $\mathbb{S}^{m}$, such that the process
$(F(t,y,z,u))_{t\in [0,T]}$ is a progressively measurable process
for each $(y,z,u)\in \mathbb{S}^{m}\times \mathbb{S}^{m}\times
\mathbb{S}^{m}$.

We consider the following matrix-valued BSDE with jump:
$$
Y_{t}=\xi+\displaystyle\int_{t}^{T}F(s,Y_{s},Z_{s},U_{s})ds-\displaystyle\int_{t}^{T}Z_{s}dW_{s}-\displaystyle\int_{t}^{T}\displaystyle\int_{E}U_{s}(e)\tilde{N}(dsde),
\eqno{(3.6)}
$$
where $\xi\in L^{2}(\Omega,\mathcal{F},P;\mathbb{S}^{m})$.

We want to study when the comparison theorem holds for two BSDE of
type (3.6).

Consider the following two BSDEs: $i=1,2,$
$$
Y^{i}_{t}=\xi^{i}+\displaystyle\int_{t}^{T}F^{i}(s,Y^{i}_{s},Z^{i}_{s},U^{i}_{s})ds-\displaystyle\int_{t}^{T}Z^{i}_{s}dW_{s}-\displaystyle\int_{t}^{T}\displaystyle\int_{E}U^{i}_{s}(e)\tilde{N}(dsde),
\eqno{(3.7)}
$$
where $F^{1},F^{2}$ satisfy (H1)-(H3), and $\xi^{1},\xi^{2}\in
L^{2}(\Omega,\mathcal{F},P;\mathbb{S}^{m})$. We study the same
problem as that in the preceding content: under which condition the
comparison theorem holds for matrix-valued BSDEs? Interestingly,
this problem is transformed again to a viability problem in
$\mathbb{S}^{m}_{+}\times\mathbb{S}^{m}$.

{\bf Theorem 3.10.} {\it Suppose that $F^{1}$ and $F^{2}$ satisfy
(H1)-(H3). Then the following are equivalent:

{\rm(i)} For any $s\in [0,T]$, $\forall \xi^{1},\xi^{2}\in
L^{2}(\Omega,\mathcal{F}_{s},P;\mathbb{S}^{m})$ such that
$\xi^{1}\geq \xi^{2}$, the unique solutions $(Y^{1},Z^{1},U^{1})$
and $(Y^{2},Z^{2},U^{2})$ in $\mathcal{B}^{2}_{[0,s]}$ valued in
$\mathbb{S}^{m}\times \mathbb{S}^{m}\times \mathbb{S}^{m}$ to the
BSDE (3.7) over time interval $[0,s]$ satisfy:
$$
Y^{1}_{t}\geq Y^{2}_{t},\mbox{ \ }t\in [0,s];
$$

{\rm(ii)} $\forall t\in [0,T]$, $\forall (y,y')$, $\forall (z,z')$,
$\forall (u,u')$,
$$
\begin{array}{ll}
&-4\langle y^{-},F^{1}(t,y^{+}+y',z,u)-F^{2}(t,y',z',u')\rangle\\
\leq & \langle
D^{2}d^{2}_{\mathbb{S}_{+}^{m}}(y)(z-z'),z-z'\rangle+C\|y^{-}\|^{2}\\
& +2\displaystyle\int_{E}
[\|(y+u(e)-u'(e))^{-}\|^{2}-\|y^{-}\|^{2}+2\langle
y^{-},u(e)-u'(e)\rangle]n(de), P-a.s.,
\end{array}
$$
where $C$ is a constant which dose not depend on
$t,(y,y'),(z,z'),(u,u')$. }

{\bf Proof: }We can see from the appendix in [10] that, for any
$y\in \mathbb{S}^{m}$, $y$ has an expression:
$$
y(\l,A)=e^{A}\sum\limits_{i=1}^{m}\l_{i}e_{i}e^{T}_{i}e^{-A},
$$
where $A$ is an antisymmetric real $m\times m$ matrix $(A^{T}=-A)$,
$\l_{i}\in R$, $\{e_{1},e_{2},...,e_{m}\}$ is the standard basis of
$R^{m}$.

If we set
$$
y^{+}(\l,A)=e^{A}\sum\limits_{i=1}^{m}\l^{+}_{i}e_{i}e^{T}_{i}e^{-A},\mbox{
\
}y^{-}(\l,A)=e^{A}\sum\limits_{i=1}^{m}\l^{-}_{i}e_{i}e^{T}_{i}e^{-A}.
$$
Then from [10], we have
$$
d^{2}_{\mathbb{S}_{+}^{m}}(y)=\|y^{-}\|^{2},
\Pi_{\mathbb{S}_{+}^{m}}(y)=y^{+},\mbox{ \ and \ }\nabla
d^{2}_{\mathbb{S}_{+}^{m}}(y)=-2y^{-},
$$
where $\|y\|=(tr(y^{2}))^{\frac{1}{2}}$. So we can use the same
method with that in Theorem 3.1 to finish the proof of this theorem.
$$
\eqno{\Box}
$$

{\bf Corollary 3.11.} {\it Let us suppose furthermore that $F^{1}$
and $F^{2}$ are independent of $z$ and $u$. Then (ii) is equivalent
to:
$$
\mbox{(ii')} -4\langle
y^{-},F^{1}(t,y^{+}+y')-F^{2}(t,y')\rangle\leq C\|y^{-}\|^{2}.
$$
}

This generalizes the result in [10].

\section{Appendix: Proofs of main theorem}\label{sec:intro}
\qquad {\bf Proof of Proposition 2.4. }Applying It\^{o}'s formula to
$e^{\b t}|Y_{t}|^{2}$, we have for $s\in [t,T]$,
$$
\begin{array}{ll}
&e^{\b s}|Y_{s}|^{2}+\displaystyle\displaystyle\int_{s}^{T}e^{\b
r}(\b
|Y_{r}|^{2}+|Z_{r}|^{2})dr+\displaystyle\int_{s}^{T}\displaystyle\int_{E}e^{\b
r}|U_{r}(e)|^{2}n(de)dr\\
=& e^{\b T}|\xi|^{2}+2\displaystyle\int_{s}^{T}e^{\b r}\langle
f(r,Y_{r},Z_{r},U_{r}),Y_{r}\rangle
dr-2\displaystyle\int_{s}^{T}e^{\b r}\langle
Y_{r},Z_{r}dW_{r}\rangle \\
&-\displaystyle\int_{s}^{T}\displaystyle\int_{E}e^{\b
r}(|Y_{r-}+U_{r}(e)|^{2}-|Y_{r-}|^{2})\tilde{N}(drde).
\end{array}
$$
Since
$$
\begin{array}{lll}
2\langle f(r,y,z,u),y\rangle &=&2\langle
f(r,y,z,u)-f(r,0,0,0),y\rangle+2\langle f(r,0,0,0),y\rangle\\
& \leq&
\frac{1}{2}|z|^{2}+\frac{1}{2}\|u\|^{2}+2(L+2L^{2})|y|^{2}+2\langle
f(r,0,0,0),y\rangle,
\end{array}
$$
and then, for $s\in [t,T]$,
$$
\begin{array}{ll}
&e^{\b s}|Y_{s}|^{2}+\displaystyle\int_{s}^{T}e^{\b
r}[(\b-2L-4L^{2})
|Y_{r}|^{2}+\frac{1}{2}|Z_{r}|^{2}]dr+\displaystyle\int_{s}^{T}\displaystyle\int_{E}\frac{1}{2}e^{\b
r}|U_{r}(e)|^{2}n(de)dr\\
\leq & e^{\b T}|\xi|^{2}+2\displaystyle\int_{s}^{T}e^{\b r}\langle
f(r,0,0,0),Y_{r}\rangle dr-2\displaystyle\int_{s}^{T}e^{\b r}\langle
Y_{r},Z_{r}dW_{r}\rangle \\
&-\displaystyle\int_{s}^{T}\displaystyle\int_{E}e^{\b
r}(|Y_{r-}+U_{r}(e)|^{2}-|Y_{r-}|^{2})\tilde{N}(drde).
\end{array}
 \eqno{(4.1)}
$$
By replacing $\b=2L+4L^{2}$, inequality (4.1) yields
$$
\begin{array}{ll}
&\displaystyle\int_{s}^{T}e^{\b
r}|Z_{r}|^{2}dr+\displaystyle\int_{s}^{T}\displaystyle\int_{E}e^{\b
r}|U_{r}(e)|^{2}n(de)dr\\
\leq & 2e^{\b T}|\xi|^{2}+4\displaystyle\int_{s}^{T}e^{\b r}\langle
f(r,0,0,0),Y_{r}\rangle dr-4\displaystyle\int_{s}^{T}e^{\b r}\langle
Y_{r},Z_{r}dW_{r}\rangle \\
&-2\displaystyle\int_{s}^{T}\displaystyle\int_{E}e^{\b
r}(|Y_{r-}+U_{r}(e)|^{2}-|Y_{r-}|^{2})\tilde{N}(drde).
\end{array}
 \eqno{(4.2)}
$$
By the Burkholder-Davis-Gundy inequality, we have
$$
\begin{array}{ll}
&E[\sup_{s\in [t,T]}|\displaystyle\int_{s}^{T}e^{\b r}\langle
Y_{r},Z_{r}dW_{r}\rangle||\mathcal{F}_{t}]\\
\leq & 2E[\sup_{s\in [t,T]}|\displaystyle\int_{t}^{s}e^{\b r}\langle
Y_{r},Z_{r}dW_{r}\rangle||\mathcal{F}_{t}]\\
\leq& C E[(\displaystyle\int^{T}_{t}e^{2\b
r}|Y_{r}|^{2}|Z_{r}|^{2}dr)^{\frac{1}{2}}|\mathcal{F}_{t}]\\
\leq & C E\{[\sup_{s\in [t,T]}(e^{\frac{1}{2}\b
s}|Y_{s}|)(\displaystyle\int_{t}^{T}e^{\b
r}|Z_{r}|^{2}dr)^{\frac{1}{2}}]|\mathcal{F}_{t}\}\\
 \leq & \frac{1}{8}E[\sup_{s\in [t,T]}e^{\b
s}|Y_{s}|^{2}|\mathcal{F}_{t}]+CE[\displaystyle\int_{t}^{T}e^{\b
r}|Z_{r}|^{2}dr|\mathcal{F}_{t}],
\end{array}
\eqno{(4.3)}
$$
here and in the sequel, $C$ denote some positive constant depending
only on $L$ and $T$, and to which we allow to change from one
formula to the other. And
$$
\begin{array}{ll}
&E[\sup_{s\in
[t,T]}\displaystyle\int_{s}^{T}\displaystyle\int_{E}e^{\b
r}(|Y_{r-}+U_{r}(e)|^{2}-|Y_{r-}|^{2})\tilde{N}(drde)|\mathcal{F}_{t}]\\
\leq& E[\displaystyle\int_{t}^{T}\displaystyle\int_{E}e^{\b
r}|U_{r}(e)|^{2}N(drde)+\displaystyle\int_{t}^{T}\displaystyle\int_{E}e^{\b
r}|U_{r}(e)|^{2}n(de)dr|\mathcal{F}_{t}]\\
&+2E\{\sup_{s\in
[t,T]}[\displaystyle\int_{t}^{s}\displaystyle\int_{E}e^{\b
r}2\langle Y_{r-},U_{r}(e)\rangle \tilde{N}(drde)|\mathcal{F}_{t}]\}\\
\leq& 2E[\displaystyle\int_{t}^{T}\displaystyle\int_{E}e^{\b
r}|U_{r}(e)|^{2}n(de)dr|\mathcal{F}_{t}]+CE[(\displaystyle\int_{t}^{T}\displaystyle\int_{E}e^{2\b
r}|Y_r|^2|U_{r}(e)|^{2}n(de)dr)^\frac{1}{2}|\mathcal{F}_{t}]\\
\leq& 2E[\displaystyle\int_{t}^{T}\displaystyle\int_{E}e^{\b
r}|U_{r}(e)|^{2}n(de)dr|\mathcal{F}_{t}]+C E[(\sup_{s\in
[t,T]}e^{\frac{1}{2}\b
s}|Y_s|)(\displaystyle\int_{t}^{T}\displaystyle\int_{E}e^{\b
r}|U_{r}(e)|^{2}n(de)dr)^\frac{1}{2}|\mathcal{F}_{t}]\\
\leq&\frac{1}{4}E[\sup_{s\in [t,T]}e^{\b
s}|Y_{s}|^{2}|\mathcal{F}_{t}]+CE[\displaystyle\int_{t}^{T}\displaystyle\int_{E}e^{\b
r}|U_{r}|^{2}n(de)dr|\mathcal{F}_{t}]
\end{array}
$$

This with (4.1), (4.2) and (4.3), we obtain
$$
\begin{array}{ll}
& \frac{1}{2}E[\sup_{s\in [t,T]}e^{\b
s}|Y_{s}|^{2}|\mathcal{F}_{t}]\\
 \leq & Ce^{\b T}E[|\xi|^{2}|\mathcal{F}_{t}]+CE[\displaystyle\int_{t}^{T}e^{\b
 r}|Y_{r}||f(r,0,0,0)|dr|\mathcal{F}_{t}]\\
\leq &Ce^{\b T}E[|\xi|^{2}|\mathcal{F}_{t}]+CE[(\sup_{s\in
[t,T]}e^{\frac{1}{2}\b
s}|Y_{s}|)\displaystyle\int_{t}^{T}e^{\frac{1}{2}\b
 r}|f(r,0,0,0)|dr|\mathcal{F}_{t}].
\end{array}
$$
Hence,
$$
E[\sup_{s\in [t,T]}e^{\b s}|Y_{s}|^{2}|\mathcal{F}_{t}]\leq Ce^{\b
T}E[|\xi|^{2}|\mathcal{F}_{t}]+CE[(\displaystyle\int_{t}^{T}e^{\frac{1}{2}\b
 r}|f(r,0,0,0)|dr)^{2}|\mathcal{F}_{t}].
$$
This, together with (4.2), gives
$$
\begin{array}{ll}
&E[\displaystyle\int_{t}^{T}e^{\b
s}|Z_{s}|^{2}ds+\displaystyle\int_{t}^{T}\displaystyle\int_{E}e^{\b
s}|U_{s}(e)|^{2}n(de)ds|\mathcal{F}_{t}]\\
\leq& Ce^{\b
T}E[|\xi|^{2}|\mathcal{F}_{t}]+CE[(\displaystyle\int_{t}^{T}e^{\frac{1}{2}\b
 r}|f(r,0,0,0)|dr)^{2}|\mathcal{F}_{t}].
\end{array}
$$
So
$$
\begin{array}{ll}
&E[\sup_{s\in [t,T]}e^{\b
s}|Y_{s}|^{2}|\mathcal{F}_{t}]+E[\displaystyle\int_{t}^{T}e^{\b
s}|Z_{s}|^{2}ds+\displaystyle\int_{t}^{T}\displaystyle\int_{E}e^{\b
s}|U_{s}(e)|^{2}n(de)ds|\mathcal{F}_{t}]\\
\leq &Ce^{\b
T}E[|\xi|^{2}|\mathcal{F}_{t}]+CE[(\displaystyle\int_{t}^{T}e^{\frac{1}{2}\b
 r}|f(r,0,0,0)|dr)^{2}|\mathcal{F}_{t}].
\end{array}
$$
Choose $C_{0}=Ce^{\b T}$, one can obtain (2.3).

In order to prove (2.4), note that BSDE (2.1) is equivalent to
$$
\bar{Y}_{t}=\displaystyle\int_{t}^{T}g(s,\bar{Y}_{s},\bar{Z}_{s},\bar{U}_{s})ds-\displaystyle\int_{t}^{T}\bar{Z}_{s}dW_{s}-\displaystyle\int_{t}^{T}\displaystyle\int_{E}\bar{U}_{s}(e)\tilde{N}(dsde),
t\in [0,T],
$$
where $\bar{Y}_{t}:=Y_{t}-E(\xi|\mathcal{F}_{t})$,
$\bar{Z}_{t}:=Z_{t}-R_{t}$, $\bar{U}_{t}:=U_{t}-V_{t}$ and
$$g(s,\bar{y},\bar{z},\bar{u}):=f(s,\bar{y}+E(\xi|\mathcal{F}_{s}),\bar{z}+R_{s},\bar{u}+V_{s}).$$
Then inequality (2.3) applied to
$(\bar{Y}_{t},\bar{Z}_{t},\bar{U}_{t})$ implies (2.4).

To prove (2.5), we set $\ep=T-t$ in (2.4) and take expectation. Then
the left thing to us is to prove, there exists a positive constant
$M_{0}$ such that
$$
E[\displaystyle\int_{t}^{T}|f(s,E(\xi|\mathcal{F}_{s}),R_{s},V_{s})|^{2}ds]\leq
M_{0}.\eqno{(4.4)}
$$
This is because
$$
\begin{array}{ll}
&E[\displaystyle\int_{t}^{T}|f(s,E(\xi|\mathcal{F}_{s}),R_{s},V_{s})|^{2}ds]\\
\leq &E[\displaystyle\int_{t}^{T}
4(|f(s,0,0,0)|^{2}+L^{2}|E(\xi|\mathcal{F}_{s})|^{2}+L^{2}|R_{s}|^{2}+L^{2}|V_{s}|^{2})ds]\\
 \leq & 4T E[\sup_{s\in
[t,T]}|f(s,0,0,0)|^{2}]+4TL^{2}E|\xi|^{2}\\
&+4L^{2}E\displaystyle\int_{t}^{T}|R_{s}|^{2}ds+4L^{2}E\displaystyle\int_{t}^{T}\displaystyle\int_{E}|V_{s}(e)|^{2}n(de)ds.
\end{array}
$$
By (H3) and Lemma 2.2, we can deduce that there exists $M_{0}>0$,
such that (4.4) holds true. So we complete the proof of the
proposition.
$$
\eqno{\Box}
$$

{\bf Proof of Theorem 2.6.} This proof is splitted into several
steps.

(a) Necessity. Let $t\in (0,T]$ and $\ep>0$ be such that, for some
$$t_{*}\geq 0,t_{\ep}:=t-\ep>t_{*}\geq 0.$$ Fix $y\in R^{m}$, $z\in
R^{m\times d}$, $u\in L^{2}(E,\mathcal{B}_{E},n; R^{m})$ and
$$
\xi:=y+z(W_{t}-W_{t_{\ep}})+\displaystyle\int_{t_{\ep}}^{t}\displaystyle\int_{E}u(e)N(drde).
$$

Denote by $(Y,Z,U)$ the unique solution to BSDE with jump
$$
Y_{s}=\xi+\displaystyle\int_{s}^{t}f(r,Y_{r},Z_{r},Z_{r})dr-\displaystyle\int_{s}^{t}Z_{r}dW_{r}-\displaystyle\int_{s}^{t}\displaystyle\int_{E}U_{r}(e)\tilde{N}(drde),s\in[t_{\ep},t].
$$

Furthermore, we introduce the Process $\hat{Y}:=\{\hat{Y}_{s},s\in
[t_{\ep},t] \}$ as follows:
$$
\begin{array}{rl}
\hat{Y}_{s}=&\xi+(t-s)f(t_{\ep},y,z,u)-z(W_{t}-W_{s})-\displaystyle\int_{s}^{t}\displaystyle\int_{E}u(e)\tilde{N}(drde)\\
=&
y+(t-s)f(t_{\ep},y,z,u)+z(W_{s}-W_{t_{\ep}})+(t-s)\displaystyle\int_{E}u(e)n(de)+\displaystyle\int^{s}_{t_{\ep}}\displaystyle\int_{E}u(e)N(drde).
\end{array}
$$

Let's compute that
$$
\begin{array}{ll}
&E[|\xi|^{2}|\mathcal{F}_{t_{*}}]\\
=& E|y+z(W_{t}-W_{t_{\ep}})+\displaystyle\int_{t_{\ep}}^{t}\displaystyle\int_{E}u(e)\tilde{N}(drde)+\displaystyle\int_{t_{\ep}}^{t}\displaystyle\int_{E}u(e)n(de)dr|^{2}\\
\leq &
4E[|y|^{2}+|z(W_{t}-W_{t_{\ep}})|^{2}+|\displaystyle\int_{t_{\ep}}^{t}\displaystyle\int_{E}u(e)\tilde{N}(drde)|^{2}+|\displaystyle\int_{t_{\ep}}^{t}\displaystyle\int_{E}u(e)n(de)dr|^{2}]\\
=& 4[|y|^{2}+\ep
z^{2}+\displaystyle\int_{t_{\ep}}^{t}\displaystyle\int_{E}|u(e)|^{2}n(de)dr+|\displaystyle\int_{t_{\ep}}^{t}\displaystyle\int_{E}u(e)n(de)dr|^{2}]\\
\leq & 4[|y|^{2}+\ep
z^{2}+(\ep+n(E)\ep^{2})\displaystyle\int_{E}|u(e)|^{2}de],
\end{array}
$$
and
$$
E[(\displaystyle\int^{t}_{t_{\ep}}|f(r,0,0,0)|dr)^{2}|\mathcal{F}_{t_{*}}]\leq
E[\ep^{2}\sup_{s\in
[t_{\ep},t]}|f(s,0,0,0)|^{2}|\mathcal{F}_{t_{*}}]
$$
So according to (2.3) in Proposition 2.4, there exists a nonnegative
random variable $\z\in L^{1}(\Omega,\mathcal{F}_{t_{*}},P)$ whose
norm depends only on $y$, $z$ and $u$, such that
$$
E[\sup_{s\in
[t_{\ep},t]}|Y_{s}|^{2}|\mathcal{F}_{t_{*}}]+E[\displaystyle\int_{t_{\ep}}^{t}|Z_{s}|^{2}ds|\mathcal{F}_{t_{*}}]
+E[\displaystyle\int_{t_{\ep}}^{t}\displaystyle\int_{E}|U_{s}(e)|^{2}n(de)ds|\mathcal{F}_{t_{*}}]\leq
\z. \eqno{(4.5)}
$$
(In the sequel, we will denote by $\z$ a nonnegative random variable
which belongs to $L^{1}(\Omega,\mathcal{F}_{t_{*}},P)$  and whose
norm depends only on $y$, $z$ and $u$, and we allow it to change
from one formula to the other and assume that $\z\geq 1, P-a.s.$ if
we need.) By (2.4) in Proposition 2.4,
$$
\begin{array}{ll}
& E[\sup_{s\in
[t_{\ep},t]}|Y_{s}-E(\xi|\mathcal{F}_{s})|^{2}|\mathcal{F}_{t_{*}}]
+E[\displaystyle\int_{t_{\ep}}^{t}|Z_{s}-z|^{2}ds|\mathcal{F}_{t_{*}}]\\
&+E[\displaystyle\int_{t_{\ep}}^{t}\displaystyle\int_{E}|U_{s}(e)-u(e)|^{2}n(de)ds|\mathcal{F}_{t_{*}}]\\
\leq & C_{0}\ep
E[\displaystyle\int_{t_{\ep}}^{t}|f(s,E(\xi|\mathcal{F}_{s}),z,u)|^{2}ds|\mathcal{F}_{t_{*}}]\\
\leq & C_{0}\ep
E[\displaystyle\int_{t_{\ep}}^{t}4(|f(s,0,0,0)|^{2}+L^{2}|E(\xi|\mathcal{F}_{s})|^{2}+L^{2}|z|^{2}+L^{2}\|u\|^{2})ds|\mathcal{F}_{t_{*}}]\\
\leq & \ep^{2}4C_{0}E[(\sup_{s\in
[t_{\ep},t]}|f(s,0,0,0)|^{2})+L^{2}|\xi|^{2}+L^{2}|z|^{2}+L^{2}\|u\|^{2}|\mathcal{F}_{t_{*}}]\\
\leq & \z\ep^{2}.
\end{array}
\eqno(4.6)
$$
 Then
$$
\begin{array}{ll}
&E[\sup_{s\in [t_{\ep},t]}|Y_{s}-y|^{2}|\mathcal{F}_{t_{*}}]\\
\leq &4E[\sup_{s\in
[t_{\ep},t]}|Y_{s}-E(\xi|\mathcal{F}_{s})|^{2}|\mathcal{F}_{t_{*}}]+4E[\sup_{s\in
[t_{\ep},t]}|z(W_{s}-W_{t_{\ep}})|^{2}|\mathcal{F}_{t_{*}}]\\
&+4E[(\ep\displaystyle\int_{E}|u(e)|n(de))^{2}|\mathcal{F}_{t_{*}}]+4E[\sup_{s\in
[t_{\ep},t]}|\displaystyle\int_{t_{\ep}}^{s}\displaystyle\int_{E}u(e)\tilde{N}(drde)|^{2}|\mathcal{F}_{t_{*}}]\\
\leq & 4[\z\ep^{2}+\ep |z|^{2}+(n(E)\ep^{2}+C\ep)\|u\|^{2}]\\
\leq & \z\ep,\mbox{ \ for \ }\ep\in (0,t-t_{*}).
\end{array}
\eqno{(4.7)}
$$

Since for $t_{\ep}\leq s\leq t$,
$$
\begin{array}{ll}
&Y_{s}-\hat{Y}_{s}\\
=&\displaystyle\int_{s}^{t}(f(r,Y_{r},Z_{r},U_{r})-f(t_{\ep},y,z,u))dr
-\displaystyle\int_{s}^{t}(Z_{r}-z)dW_{r}-\displaystyle\int_{s}^{t}\displaystyle\int_{E}(U_{r}(e)-u(e))\tilde{N}(drde),
\end{array}
$$
we can apply It\^{o}'s formula to $|Y_{s}-\hat{Y}_{s}|^{2}$ over
$[t_{\ep},t]$ and with the same technique as that in Proposition
2.4, noting (4.6) and (4.7), we have
$$
\begin{array}{ll}
&
E[\sup_{s\in[t_{\ep},t]}|Y_{s}-\hat{Y}_{s}|^{2}|\mathcal{F}_{t_{*}}]+E[\displaystyle\int^{t}_{t_{\ep}}|Z_{r}-z|^{2}dr|\mathcal{F}_{t_{*}}]\\
&+E[\displaystyle\int^{t}_{t_{\ep}}\displaystyle\int_{E}|U_{r}(e)-u(e)|^{2}n(de)dr|\mathcal{F}_{t_{*}}]\\
\leq & C\ep
E[\displaystyle\int^{t}_{t_{\ep}}|f(r,Y_{r},Z_{r},U_{r})-f(t_{\ep},y,z,u)|^{2}dr|\mathcal{F}_{t_{*}}]\\
\leq &
4C\ep^{2}\{E[\sup_{s\in[t_{\ep},t]}|f(s,y,z,u)-f(t_{\ep},y,z,u)|^{2}|\mathcal{F}_{t_{*}}]
+L^{2}E[\sup_{s\in[t_{\ep},t]}|Y_{s}-y|^{2}|\mathcal{F}_{t_{*}}]\}\\
&+4CL^{2}\ep\{E[\displaystyle\int^{t}_{t_{\ep}}|Z_{r}-z|^{2}dr|\mathcal{F}_{t_{*}}]+E[\displaystyle\int^{t}_{t_{\ep}}\displaystyle\int_{E}|U_{r}(e)-u(e)|^{2}n(de)dr|\mathcal{F}_{t_{*}}]\}\\
\leq & \z\ep^{2}\b^{1}_{\ep},
\end{array}
\eqno{(4.8)}
$$
where
$$
\b^{1}_{\ep}=8CL^{2}\ep+4\frac{C}{\z}E[\sup_{s\in[t_{\ep},t]}|f(s,y,z,u)-f(t_{\ep},y,z,u)|^{2}|\mathcal{F}_{t_{*}}].
$$

Observe that by (H1):
$$
\sup_{s\in[t_{\ep},t]}|f(s,y,z,u)-f(t_{\ep},y,z,u)|^{2}\rightarrow
0,\mbox{ \ when \ }\ep\rightarrow 0, P-a.s.,
$$
and that
$$
\begin{array}{ll}
&\sup_{s\in[t_{\ep},t]}|f(s,y,z,u)-f(t_{\ep},y,z,u)|^{2}\\
\leq & 4\sup _{t\in [0,T]}|f(t,y,z,u)|^{2}\\
\leq & 8\sup_{t\in
[0,T]}|f(t,0,0,0)|^{2}+8L^{2}(|y|^{2}+|z|^{2}+\|u\|^{2})\in
L^{1}(\Omega,\mathcal{F}_{T},P).
\end{array}
$$
Hence from Lebesgue's dominated convergence theorem,
$\b^{1}_{\ep}\rightarrow 0, P-a.s.$, as $\ep\rightarrow 0$.

Let us now establish two auxiliary results on the processes $Y$ and
$\hat{Y}$ which will enable us to finish the proof of the necessity.

{\bf Lemma 4.1.} {\it Under the assumptions made above, there is
some nonnegative random variable $\z\in
L^{1}(\Omega,\mathcal{F}_{t_{*}},P)$ whose norm depends only on $y$,
$z$ and $u$ such that
$$
E[|d^{2}_{K}(Y_{t_{\ep}})-d^{2}_{K}(\hat{Y}_{t_{\ep}})||\mathcal{F}_{t_{*}}]\leq
\z\ep\sqrt{\b^{1}_{\ep}},
$$
for any $\ep>0$ with $t-\ep>t_{*}$. }

{\bf Proof of Lemma 4.1.} Note that, since $K\neq {\O}$, there
exists some constant $C>0$ such that, $\forall x,x'\in R^{m}$
$$
|d^{2}_{K}(x)-d^{2}_{K}(x')|=|d_{K}(x)-d_{K}(x')|(d_{K}(x)+d_{K}(x'))
\leq C(1+|x|+|x'|)|x-x'|.
$$
So from (4.8)
$$
\begin{array}{ll}
&E[|d^{2}_{K}(Y_{t_{\ep}})-d^{2}_{K}(\hat{Y}_{t_{\ep}})||\mathcal{F}_{t_{*}}]\\
\leq &
CE[(1+|Y_{t_{\ep}}|+|\hat{Y}_{t_{\ep}}|)|Y_{t_{\ep}}-\hat{Y}_{t_{\ep}}||\mathcal{F}_{t_{*}}]\\
\leq &
C(E[(1+|Y_{t_{\ep}}|+|\hat{Y}_{t_{\ep}}|)^{2}|\mathcal{F}_{t_{*}}])^{\frac{1}{2}}(E[|Y_{t_{\ep}}-\hat{Y}_{t_{\ep}}|^{2}|\mathcal{F}_{t_{*}}])^{\frac{1}{2}}\\
\leq &
C\ep\sqrt{\z\b^{1}_{\ep}}(E[(1+|Y_{t_{\ep}}|+|\hat{Y}_{t_{\ep}}|)^{2}|\mathcal{F}_{t_{*}}])^{\frac{1}{2}}.
\end{array}
$$
Observe that $\hat{Y}_{t_{\ep}}=y+\ep
f(t_{\ep},y,z,u)+\ep\displaystyle\int_{E}u(e)n(de)$, noting (H3) and
(4.5), we have
$$
E[|d^{2}_{K}(Y_{t_{\ep}})-d^{2}_{K}(\hat{Y}_{t_{\ep}})||\mathcal{F}_{t_{*}}]\leq
\z\ep\sqrt{\b^{1}_{\ep}},
$$
for some $\z\in L^{1}(\Omega,\mathcal{F}_{t_{*}},P)$. This completes
the proof of Lemma 4.1.
$$
\eqno{\Box}$$

{\bf Lemma 4.2.} {\it We can find some
$\mathcal{F}_{t_{*}}$-measurable random variable
$\g_{\ep}=\g_{\ep}(y,z,u)$ with $\liminf_{\ep\rightarrow
0}\g_{\ep}\geq 0$ such that
$$
\begin{array}{ll}
&E[d^{2}_{K}(\hat{Y}_{t_{\ep}})-d^{2}_{K}(\xi)|\mathcal{F}_{t_{*}}]\\
\geq & \ep\{E[\langle\nabla
d^{2}_{K}(y),f(t,y,z,u)+\displaystyle\int_{E}u(e)n(de)\rangle|\mathcal{F}_{t_{*}}]\\
& \mbox{ \ \ }-\frac{1}{2}\langle
D^{2}[d^{2}_{K}(y)]z,z\rangle-\frac{1}{\ep}E[\displaystyle\int_{t_{\ep}}^{t}\displaystyle\int_{E}(d^{2}_{K}(\xi_{s-}+u(e))-d^{2}_{K}(\xi_{s-}))n(de)ds|\mathcal{F}_{t_{*}}]+\g_{\ep}\},
\end{array}
$$
where
$\xi_{s}=y+z(W_{s}-W_{t_{\ep}})+\displaystyle\int_{t_{\ep}}^{s}\displaystyle\int_{E}u(e)N(drde),s\in
[t_{\ep},t]$. }

{\bf Proof: }We know that the function $d^{2}_{K}(\cdot)$ is twice
differentiable almost everywhere. Let us denote by $\Lambda_{K}$ the
set of all points of $R^{m}$ where $d^{2}_{K}$ is twice
differentiable. This set is of full Lebesgue measure. Let us fix now
$y\in \Lambda_{K}$ and define the following function $\a :
R^{m}\rightarrow R$:
$$
\a(x):=d^{2}_{K}(x+y)-d^{2}_{K}(y)-\langle \nabla
d^{2}_{K}(y),x\rangle-\frac{1}{2}\langle
D^{2}[d^{2}_{K}(y)]x,x\rangle.
$$
There are two following properties of $\a(\cdot)$ (see [7]):
$$
\lim\limits_{|x|\rightarrow 0}\frac{\a(x)}{|x|^{2}}=0,
$$

$$
\forall x\in R^{m},\mbox{ \ }\a(x)\leq
|x|^{2}(1+|D^{2}d^{2}_{K}(y)|). \eqno{(4.9)}
$$

First we substitute $$x=\hat{Y}_{t_{\ep}}-y=\ep
f(t_{\ep},y,z,u)+\ep\displaystyle\int_{E}u(e)n(de)$$ in the
definition of $\a(\cdot)$. This provides us
$$
E[d^{2}_{K}(\hat{Y}_{t_{\ep}})|\mathcal{F}_{t_{*}}]=d^{2}_{K}(y)+\ep
E[\langle\nabla
d^{2}_{K}(y),f(t_{\ep},y,z,u)+\displaystyle\int_{E}u(e)n(de)\rangle|\mathcal{F}_{t_{*}}]+\ep\g^{1}_{\ep},
$$
where
$$
\begin{array}{l}
\g^{1}_{\ep}:=\frac{\ep}{2}E[\langle
D^{2}d^{2}_{K}(y)(f(t_{\ep},y,z,u)+\displaystyle\int_{E}u(e)n(de)),f(t_{\ep},y,z,u)+\displaystyle\int_{E}u(e)n(de)
\rangle\\
\mbox{ \ \ \ \ \ \ \ \ \ \ \ \ }+\frac{2}{\ep^{2}}\a(\ep
f(t_{\ep},y,z,u)+\ep\displaystyle\int_{E}u(e)n(de))|\mathcal{F}_{t_{*}}].
\end{array}
$$
By (H3), $f(t_{\ep},y,z,u)+\displaystyle\int_{E}u(e)n(de)$ can be
dominated by some nonnegative random variable which belongs to
$L^{2}(\Omega,\mathcal{F}_{t_{*}},P)$, so there exists some
nonnegative random variable $\z\in
L^{1}(\Omega,\mathcal{F}_{t_{*}},P)$ whose norm depends only on $y$,
$z$ and $u$ such that
$$
|\g^{1}_{\ep}|\leq \z\ep+E[\frac{\a(
x)}{|x|^{2}}\ep|f(t_{\ep},y,z,u)+\displaystyle\int_{E}u(e)n(de)|^{2}|\mathcal{F}_{t_{*}}].
$$
Hence from Lebesgue's dominated convergence theorem,
$\g^{1}_{\ep}\rightarrow 0,P-a.s.$, as $\ep$ tends to 0.

We now substitute
$$x=\xi-y=z(W_{t}-W_{t_{\ep}})+\displaystyle\int_{t_{\ep}}^{t}\displaystyle\int_{E}u(e)N(drde)$$
in the definition of $\a(\cdot)$. If we set
$$\xi_{s}=y+z(W_{s}-W_{t_{\ep}})+\displaystyle\int_{t_{\ep}}^{s}\displaystyle\int_{E}u(e)N(drde),s\in
[t_{\ep},t],$$ we have
$$
\begin{array}{ll}
&E[d^{2}_{K}(\xi)-d^{2}_{K}(y)|\mathcal{F}_{t_{*}}]\\
=&
E\{\displaystyle\int^{t}_{t_{\ep}}\displaystyle\int_{E}[d^{2}_{K}(\xi_{s-}+u(e))-d^{2}_{K}(\xi_{s-})]N(dsde)\\
&\mbox{ \ \ \ }+\frac{1}{2}\langle
D^{2}d^{2}_{K}(y)z(W_{t}-W_{t_{\ep}}),z(W_{t}-W_{t_{\ep}})\rangle+\a(z(W_{t}-W_{t_{\ep}}))|\mathcal{F}_{t_{*}}\}\\
=&\ep
E\{\frac{1}{\ep}\displaystyle\int^{t}_{t_{\ep}}\displaystyle\int_{E}[d^{2}_{K}(\xi_{s-}+u(e))-d^{2}_{K}(\xi_{s-})]n(de)ds+\frac{1}{2}\langle
D^{2}d^{2}_{K}(y)z,z\rangle|\mathcal{F}_{t_{*}}\}+\ep\g^{2}_{\ep},
\end{array}
$$
where
$$
\g^{2}_{\ep}=\frac{1}{\ep}E[\a(z(W_{t}-W_{t_{\ep}}))|\mathcal{F}_{t_{*}}]=\frac{1}{\ep}E[\a(\sqrt{\ep}zW_{1})]
$$
is such that
$$
\limsup_{\ep\rightarrow 0}\g^{2}_{\ep}\leq 0.
$$
In fact, on one hand,
$$
\frac{1}{\ep}\a(\sqrt{\ep}zW_{1})\rightarrow 0, P-a.s., \mbox{ \ as
\ }\ep\rightarrow 0.
$$
And on the other hand, for $\ep>0$, with (4.9),
$$
\frac{1}{\ep}\a(\sqrt{\ep}zW_{1})\leq
(1+|D^{2}d^{2}_{K}(y)|)|z|^{2}|W_{1}|^{2}\in L^{1}(P).
$$
Finally, we get from Fatou's Lemma
$$
\limsup_{\ep\rightarrow 0}\g^{2}_{\ep}\leq E[\limsup_{\ep\rightarrow
0}\frac{1}{\ep}\a(\sqrt{\ep}zW_{1})|\mathcal{F}_{t_{*}}]=0,P-a.s..
$$
Therefore
$$
\begin{array}{ll}
&E[d^{2}_{K}(\hat{Y}_{t_{\ep}})-d^{2}_{K}(\xi)|\mathcal{F}_{t_{*}}]\\
=& \ep\{E[\langle\nabla
d^{2}_{K}(y),f(t,y,z,u)+\displaystyle\int_{E}u(e)n(de)\rangle|\mathcal{F}_{t_{*}}]\\
& \mbox{ \ \ }+E[\langle\nabla
d^{2}_{K}(y),f(t_{\ep},y,z,u)-f(t,y,z,u)\rangle|\mathcal{F}_{t_{*}}]-\frac{1}{2}\langle
D^{2}[d^{2}_{K}(y)]z,z\rangle\\
&\mbox{ \ \ }-\frac{1}{\ep}E[\displaystyle\int_{t_{\ep}}^{t}\displaystyle\int_{E}(d^{2}_{K}(\xi_{s-}+u(e))-d^{2}_{K}(\xi_{s-}))n(de)ds|\mathcal{F}_{t_{*}}]+\g^{1}_{\ep}-\g^{2}_{\ep}\}\\
\geq & \ep\{E[\langle\nabla
d^{2}_{K}(y),f(t,y,z,u)+\displaystyle\int_{E}u(e)n(de)\rangle|\mathcal{F}_{t_{*}}]-\frac{1}{2}\langle
D^{2}[d^{2}_{K}(y)]z,z\rangle\\
& \mbox{ \ \
}-\frac{1}{\ep}E[\displaystyle\int_{t_{\ep}}^{t}\displaystyle\int_{E}(d^{2}_{K}(\xi_{s-}+u(e))-d^{2}_{K}(\xi_{s-}))n(de)ds|\mathcal{F}_{t_{*}}]+\g^{1}_{\ep}-\g^{2}_{\ep}-\g^{3}_{\ep}\},
\end{array}
$$
where $\g^{3}_{\ep}=E[|\nabla
d^{2}_{K}(y)||f(t_{\ep},y,z,u)-f(t,y,z,u)||\mathcal{F}_{t_{*}}]$. By
(H1) and (H3), we have
$$
\lim_{\ep\rightarrow 0}\g^{3}_{\ep}\rightarrow 0,P-a.s..
$$
So the proof of Lemma 4.2 is completed by setting
$\g_{\ep}=\g^{1}_{\ep}-\g^{2}_{\ep}-\g^{3}_{\ep}$.
$$
\eqno{\Box}
$$

Note that due to the Lemma 4.1 and 4.2:
$$
\begin{array}{ll}
&E[d^{2}_{K}(Y_{t_{\ep}})-d^{2}_{K}(\xi)|\mathcal{F}_{t_{*}}]\\
\geq & \ep\{E[\langle\nabla
d^{2}_{K}(y),f(t,y,z,u)+\displaystyle\int_{E}u(e)n(de)\rangle|\mathcal{F}_{t_{*}}]\\
& \mbox{ \ \ }-\frac{1}{2}\langle
D^{2}[d^{2}_{K}(y)]z,z\rangle-\frac{1}{\ep}E[\displaystyle\int_{t_{\ep}}^{t}\displaystyle\int_{E}(d^{2}_{K}(\xi_{s-}+u(e))-d^{2}_{K}(\xi_{s-}))n(de)ds|\mathcal{F}_{t_{*}}]+\g_{\ep}\},
\end{array}
$$
for some $\g_{\ep}=\g_{\ep}(y,z,u)$ such that
$\liminf_{\ep\rightarrow 0}\g_{\ep}\geq 0$, $P-a.s.$.

Let us return to the proof of the necessity. For this denote by
$(\tilde{Y},\tilde{Z},\tilde{U})$ the unique solution to the
following BSDE with jump:
$$
\tilde{Y}_{s}=\eta+\displaystyle\int_{s}^{t}f(r,\tilde{Y}_{r},\tilde{Z}_{r},\tilde{U}_{r})dr-\displaystyle\int_{s}^{t}\tilde{Z}_{r}dW_{r}-\displaystyle\int_{s}^{t}\displaystyle\int_{E}\tilde{U}_{r}(e)\tilde{N}(drde),s\in
[t_{\ep},t],
$$
where $\eta\in L^{2}(\Omega,\mathcal{F},P)$ is a measurable
selection of the set
$$
\{(\omega,x)\in \Omega\times R^{m}|x\in \Pi_{K}(\xi(\omega))\}\in
\mathcal{F}_{t}\otimes \mathcal{B}(R^{m}).
$$

We assume that $K$ enjoys the BSVP. Hence $\tilde{Y}_{s}\in K$, for
$t_{\ep}\leq s\leq t$, $P-a.s.$. This implies
$$
0 \geq
E[d^{2}_{K}(Y_{t_{\ep}})-|Y_{t_{\ep}}-\tilde{Y}_{t_{\ep}}|^{2}|\mathcal{F}_{t_{*}}]=
E[d^{2}_{K}(Y_{t_{\ep}})-d^{2}_{K}(\xi)|\mathcal{F}_{t_{*}}]-E[|Y_{t_{\ep}}-\tilde{Y}_{t_{\ep}}|^{2}-|\xi-\eta|^{2}|\mathcal{F}_{t_{*}}].
$$
From It\^{o}'s formula,
$$
\begin{array}{ll}
&E[|Y_{t_{\ep}}-\tilde{Y}_{t_{\ep}}|^{2}-|\xi-\eta|^{2}|\mathcal{F}_{t_{*}}]\\
=& 2E[\displaystyle\int_{t_{\ep}}^{t}\langle
Y_{s}-\tilde{Y}_{s},f(s,Y_{s},Z_{s},U_{s})-f(s,\tilde{Y}_{s},\tilde{Z}_{s},\tilde{U}_{s})\rangle
ds|\mathcal{F}_{t_{*}}]\\
&\mbox{ \ \ \
}-E[\displaystyle\int^{t}_{t_{e}}|Z_{s}-\tilde{Z}_{s}|^{2}ds+\displaystyle\int^{t}_{t_{e}}\displaystyle\int_{E}|U_{s}-\tilde{U}_{s}|^{2}n(de)ds|\mathcal{F}_{t_{*}}]\\
\leq & C
\displaystyle\int^{t}_{t_{e}}E[|Y_{s}-\tilde{Y}_{s}|^{2}|\mathcal{F}_{t_{*}}]ds,
\end{array}
$$
where $C=2L+2L^{2}$. Consequently
$$
\begin{array}{ll}
0& \geq
\frac{1}{\ep}E[d^{2}_{K}(Y_{t_{\ep}})-d^{2}_{K}(\xi)|\mathcal{F}_{t_{*}}]-\frac{C}{\ep}\displaystyle\int_{t_{\ep}}^{t}E[|Y_{s}-\tilde{Y}_{s}|^{2}|\mathcal{F}_{t_{*}}]ds\\
& \geq E[\langle\nabla
d^{2}_{K}(y),f(t,y,z,u)+\displaystyle\int_{E}u(e)n(de)\rangle|\mathcal{F}_{t_{*}}]
-\frac{1}{2}\langle
D^{2}[d^{2}_{K}(y)]z,z\rangle+\g_{\ep}\\
& \mbox{ \ \ }
-\frac{1}{\ep}E[\displaystyle\int_{t_{\ep}}^{t}\displaystyle\int_{E}(d^{2}_{K}(\xi_{s-}+u(e))-d^{2}_{K}(\xi_{s-}))n(de)ds|\mathcal{F}_{t_{*}}]
-\frac{C}{\ep}\displaystyle\int_{t_{\ep}}^{t}E[|Y_{s}-\tilde{Y}_{s}|^{2}|\mathcal{F}_{t_{*}}]ds.
\end{array}
$$
In order to estimate the integral term in the last estimate, note
that
$$
\begin{array}{ll}
& E[|Y_{s}-\tilde{Y}_{s}|^{2}|\mathcal{F}_{t_{*}}]\\
=&
E[|(Y_{s}-\xi)-(\tilde{Y}_{s}-\eta)+(\xi-\eta)|^{2}|\mathcal{F}_{t_{*}}]\\
\leq &
E[d^{2}_{K}(\xi)|\mathcal{F}_{t_{*}}]+2E[|Y_{s}-\xi|^{2}+|\tilde{Y}_{s}-\eta|^{2}|\mathcal{F}_{t_{*}}]\\
&+2(E[d^{2}_{K}(\xi)|\mathcal{F}_{t_{*}}])^{\frac{1}{2}}\{(E[|Y_{s}-\xi|^{2}|\mathcal{F}_{t_{*}}])^{\frac{1}{2}}
+(E[|\tilde{Y}_{s}-\eta|^{2}|\mathcal{F}_{t_{*}}])^{\frac{1}{2}}\}.
\end{array}
\eqno{(4.10)}
$$
While
$$
\begin{array}{ll}
&|Y_{s}-\xi|^{2}+|\tilde{Y}_{s}-\eta|^{2}\\
\leq &
2(|Y_{s}-E[\xi|\mathcal{F}_{s}]|^{2}+|\xi-E[\xi|\mathcal{F}_{s}]|^{2}+|\tilde{Y}_{s}-E[\eta|\mathcal{F}_{s}]|^{2}+|\eta-E[\eta|\mathcal{F}_{s}]|^{2}),
\end{array}
$$
so from Proposition 2.4 we conclude
$$
\sup_{s\in
[t_{\ep},t]}E[|Y_{s}-\xi|^{2}+|\tilde{Y}_{s}-\eta|^{2}|\mathcal{F}_{t_{*}}]\rightarrow
0,P-a.s.,
$$
as $\ep$ tends to 0. On the other hand, $$\begin{array}{ll}
&E[d^{2}_{K}(\xi)-d^{2}_{K}(y)|\mathcal{F}_{t_{*}}]\\
 =&
E\{\displaystyle\int^{t}_{t_{\ep}}\displaystyle\int_{E}[d^{2}_{K}(\xi_{s-}+u(e))-d^{2}_{K}(\xi_{s-})]N(dsde)+d^{2}_{K}(y+z(W_{t}-W_{t_{\ep}}))-d^{2}_{K}(y)|\mathcal{F}_{t_{*}}\}.
\end{array}$$ Recall that
$\xi_{s}=y+z(W_{s}-W_{t_{\ep}})+\displaystyle\int_{t_{\ep}}^{s}\displaystyle\int_{E}u(e)N(drde),s\in
[t_{\ep},t]$, so
$$
\xi_{s-}\rightarrow y,P-a.s.,\mbox{ \ as \ }\ep\rightarrow 0.
$$
By the Lebesgue theorem of dominated convergence, when $\ep$ tends
to 0,
$$
\begin{array}{ll}
&\frac{1}{\ep}E\{\displaystyle\int^{t}_{t_{\ep}}\displaystyle\int_{E}[d^{2}_{K}(\xi_{s-}+u(e))-d^{2}_{K}(\xi_{s-})]N(dsde)|\mathcal{F}_{t_{*}}\}\\
\rightarrow&
\displaystyle\int_{E}[d^{2}_{K}(y+u(e))-d^{2}_{K}(y)]n(de),P-a.s..
\end{array}\eqno{(4.11)}
$$
Thus
$$
E\{\displaystyle\int^{t}_{t_{\ep}}\displaystyle\int_{E}[d^{2}_{K}(\xi_{s-}+u(e))-d^{2}_{K}(\xi_{s-})]N(dsde)|\mathcal{F}_{t_{*}}\}\rightarrow
0,P-a.s., \mbox{ \ as \ }\ep\rightarrow 0.
$$
Obviously, the function $d^{2}_{K}$ is continuous at $y$, using
again Lebesgue's dominated convergence theorem, we have
$$
E[d^{2}_{K}(y+z(W_{t}-W_{t_{\ep}}))-d^{2}_{K}(y)|\mathcal{F}_{t_{*}}]\rightarrow
0, P-a.s.,\mbox{ \ as \ }\ep\rightarrow 0.
$$
So we conclude that
$$
E[d^{2}_{K}(\xi)-d^{2}_{K}(y)|\mathcal{F}_{t_{*}}]\rightarrow 0,
P-a.s.,\mbox{ \ as \ }\ep\rightarrow 0.
$$

Consequently, from (4.10) and above, we have
$$
E[|Y_{s}-\tilde{Y}_{s}|^{2}|\mathcal{F}_{t_{*}}]\leq
d^{2}_{K}(y)+\b^{2}_{\ep},s\in [t_{\ep},t],
$$
where $\b^{2}_{\ep}$ converges to 0, $P-a.s.$, as $\ep$ tends to 0.

Therefore
$$
\frac{C}{\ep}\displaystyle\int_{t_{\ep}}^{t}E[|Y_{s}-\tilde{Y}_{s}|^{2}|\mathcal{F}_{t_{*}}]ds\leq
C(d^{2}_{K}(y)+\b^{2}_{\ep})
$$
and
$$
\begin{array}{l} E[\langle\nabla
d^{2}_{K}(y),f(t,y,z,u)+\displaystyle\int_{E}u(e)n(de)\rangle|\mathcal{F}_{t_{*}}]
-\frac{1}{2}\langle
D^{2}[d^{2}_{K}(y)]z,z\rangle\\
-\frac{1}{\ep}E[\displaystyle\int_{t_{\ep}}^{t}\displaystyle\int_{E}(d^{2}_{K}(\xi_{s-}+u(e))-d^{2}_{K}(\xi_{s-}))n(de)ds|\mathcal{F}_{t_{*}}]
-Cd^{2}_{K}(y)+\g_{\ep}-C\b^{2}_{\ep}\leq 0,P-a.s..
\end{array}
$$
Finally, since $\liminf_{\ep\rightarrow
0}(\g_{\ep}-C\b^{2}_{\ep})\geq 0, P-a.s.$. This with (4.11), let
$\ep$ tend to 0, we have
$$
\begin{array}{l} E[\langle\nabla
d^{2}_{K}(y),f(t,y,z,u)+\displaystyle\int_{E}u(e)n(de)\rangle|\mathcal{F}_{t_{*}}]
\\
-\frac{1}{2}\langle
D^{2}[d^{2}_{K}(y)]z,z\rangle-\displaystyle\int_{E}[d^{2}_{K}(y+u(e))-d^{2}_{K}(y)]n(de)
-Cd^{2}_{K}(y)\leq 0,
\end{array}
$$
$P$-almost everywhere, for $t_{*}\in [0,T)$. Passing to the limit
$t_{*}\rightarrow t$, we obtain the wished result if we choose
$C^{*}=2C$ and note that
$$
\nabla d^{2}_{K}(y)=2(y-\Pi_{K}(y)).
$$
$$
\eqno{\Box}
$$
(b)Sufficiency. Let $K$ be a nonempty convex closed subset of
$R^{m}$. Suppose that (2.6) holds true. Let $\eta\in
C^{\infty}(R^{m})$ be a nonnegative function with support in the
unit ball and such that
$$
\displaystyle\int_{R^{m}}\eta(x)dx=1.
$$

For $\d>0$, we put
$$
\begin{array}{l}
\eta_{\d}(x):=\frac{1}{\d^{m}}\eta(\frac{x}{\d})\\
\ph_{\d}(x):=d^{2}_{K}\star
\eta_{\d}(x):=\displaystyle\int_{R^{m}}d^{2}_{K}(x-x')\eta_{\d}(x')dx'=\displaystyle\int_{R^{m}}d^{2}_{K}(x')\eta_{\d}(x-x')dx',x\in
R^{m}.
\end{array}
$$
Obviously, $\ph_{\d}\in C^{\infty}(R^{m})$. We can see from [7] that
the function $\ph_{\d}$ satisfies the following properties
$$
\left\{ \begin{array}{ll} {\rm(i)}& 0\leq \ph_{\d}(x)\leq
(d_{K}(x)+\d)^{2},\\
{\rm(ii)}& \nabla \ph_{\d}(x)=\displaystyle\int_{R^{m}}(\nabla
d^{2}_{K})(x')\eta_{\d}(x-x')dx',\\
& |\nabla \ph_{\d}(x)|\leq 2(d_{K}(x)+\d),\\
{\rm(iii)}&
D^{2}\ph_{\d}(x)=\displaystyle\int_{R^{m}}D^{2}[d^{2}_{K}](x')\eta_{\d}(x-x')dx',\\
& 0\leq D^{2}\ph_{\d}(x)\leq 2I,
\end{array}
\right. \eqno{(4.12)}$$ for all $x\in R^{m}$.

Consider $\xi\in L^{2}(\Omega,\mathcal{F}_{T},P,K)$ and let
$(Y,Z,U)$ be the unique solution to the following BSDE with jump
$$
Y_{t}=\xi+\displaystyle\int_{t}^{T}f(s,Y_{s},Z_{s},U_{s})ds-\displaystyle\int_{t}^{T}Z_{s}dW_{s}
-\displaystyle\int_{t}^{T}\displaystyle\int_{E}U_s(e)\tilde{N}(dsde),\mbox{
\ } t\in [0,T].
$$
(4.12) enables us to apply It\^{o}'s formula to $\ph_{\d}(Y_{t})$
and to deduce that, for $0\leq t\leq T$, $\d>0$,
$$
\begin{array}{ll}
&E\ph_{\d}(Y_{t})\\
= &E\ph_{\d}(\xi)+ E\displaystyle\int_{t}^{T}\langle (\nabla
\ph_{\d})(Y_{s}),f(s,Y_{s},Z_{s},U_{s})\rangle
ds-\frac{1}{2}E\displaystyle\int_{t}^{T}\langle
(D^{2}\ph_{\d})(Y_{s})Z_{s},Z_{s}\rangle
ds\\
&-E\displaystyle\int_{t}^{T}\displaystyle\int_{E}[\ph_{\d}(Y_{s}+U_{s}(e))-\ph_{\d}(Y_{s})-\langle
\nabla \ph_{\d}(Y_{s}),U_{s}(e)\rangle]n(de)ds\\
\leq &
\d^{2}-E\displaystyle\int_{t}^{T}\displaystyle\int_{R^{m}}\{\langle
\nabla
d^{2}_{K}(y),f(s,y,Z_{s},U_{s})-f(s,Y_{s},Z_{s},U_{s}))\rangle\eta_{\d}(Y_{s}-y)\}dyds\\
& +E\displaystyle\int_{t}^{T}\displaystyle\int_{R^{m}} \{\langle
\nabla d^{2}_{K}(y),f(s,y,Z_{s},U_{s})\rangle-\frac{1}{2}\langle
D^{2}(d^{2}_{K}(y))Z_{s},Z_{s}\rangle\\
&\mbox{ \ \ \ \ \ \ \ \ \ \ \ \ \ \
}-\displaystyle\int_{E}[d^{2}_{K}(y+U_{s}(e))-d^{2}_{K}(y)-\langle\nabla
d^{2}_{K}(y),U_{s}(e)\rangle]n(de)\}\eta_{\d}(Y_{s}-y)dyds.
\end{array}
$$
Then from (2.6) and (4.12), for $\d\in (0,1)$,
$$
\begin{array}{ll}
&E\ph_{\d}(Y_{t})\\
\leq &\d^{2}+
\frac{C^{*}}{2}E\displaystyle\int_{t}^{T}\displaystyle\int_{R^{m}}d^{2}_{K}(y)\eta_{\d}(Y_{s}-y)dyds\\
&+E\displaystyle\int_{t}^{T}\displaystyle\int_{R^{m}}2d_{K}(y)\max\limits_{y:|y-Y_{s}|\leq
\d}|f(s,y,Z_{s},U_{s})-f(s,Y_{s},Z_{s},U_{s}))|\eta_{\d}(Y_{s}-y)dyds\\
\leq &\d^{2}+
\frac{C^{*}}{2}\displaystyle\int_{t}^{T}E[\ph_{\d}(Y_{s})]ds\\
&+E\displaystyle\int_{t}^{T}(1+\ph_{\d}(Y_{s}))\max\limits_{y:|y-Y_{s}|\leq
\d}|f(s,y,Z_{s},U_{s})-f(s,Y_{s},Z_{s},U_{s}))|ds.
\end{array}
$$
Since
$$
\begin{array}{ll}
&
E\displaystyle\int_{t}^{T}\ph_{\d}(Y_{s})\max\limits_{y:|y-Y_{s}|\leq
\d}|f(s,y,Z_{s},U_{s})-f(s,Y_{s},Z_{s},U_{s}))|ds\\
\leq &
E\displaystyle\int_{t}^{T}\ph_{\d}(Y_{s})\max\limits_{y:|y-Y_{s}|\leq
\d}L|Y_{s}-y|ds\\
\leq &E\displaystyle\int_{t}^{T}L\d\ph_{\d}(Y_{s})ds,
\end{array}
$$ We choose $\d< \frac{1}{L}$, then
$$
E\ph_{\d}(Y_{t})\leq
\d^{2}+L(T-t)\d+(1+\frac{C^{*}}{2})\displaystyle\int_{t}^{T}E\ph_{\d}(Y_{s})ds.
\eqno{(4.13)}
$$
for $0\leq t\leq T$. On the other hand, from (4.12-i), we deduce
$$
E\ph_{\d}(Y_{t})\leq E[d_{K}(Y_{t})+\d]^{2}\leq
2E(d^{2}_{K}(Y_{t})+\d^{2})\leq
2E(2|Y_{t}|^{2}+2d^{2}_{K}(0)+\d^{2})<+\infty.
$$
This allows to apply Gronwall's inequality to (4.13), it yields
$$
E\ph_{\d}(Y_{t})\leq [\d^{2}+L(T-t)\d]e^{(1+\frac{C^{*}}{2})T}.
$$
Note that
$$
\ph_{\d}(Y_{t})=\displaystyle\int_{R^{m}}d^{2}_{K}(Y_{t}-y)\eta_{\d}(y)dy\rightarrow
d^{2}_{K}(Y_{t}),\mbox{ \ as \ }\d\rightarrow 0, P-a.s..
$$
Thus, from Fatou's Lemma, we conclude that
$$
Ed^{2}_{K}(Y_{t})=E[\lim\limits_{\d\rightarrow
0}\ph_{\d}(Y_{t})]\leq \liminf\limits_{\d \rightarrow
0}E\ph_{\d}(Y_{t})\leq \liminf\limits_{\d \rightarrow
0}[\d^{2}+L(T-t)\d]e^{(1+\frac{C^{*}}{2})T}=0,
$$
for $0\leq t\leq T$. That is
$$
P\{\omega: Y_{t}(\omega)\in K\}=1, \mbox{ \ }\forall t\in [0,T]
$$
which is equivalent to
$$
P\{\omega: Y_{t}(\omega)\in K, \forall t\in [0,T]\}=1.
$$
$$
\eqno{\Box}
$$

\end{document}